\documentclass[11pt,a4paper]{article}

\usepackage[a4paper,margin=1in]{geometry}
\usepackage[T1]{fontenc}
\usepackage[utf8]{inputenc}
\usepackage{lmodern}
\usepackage[english]{babel}
\usepackage{graphicx}
\usepackage{enumitem}
\usepackage{mathrsfs}
\usepackage{amsfonts,amsmath,amssymb,amsthm,bm}
\usepackage{xcolor}
\usepackage{mathtools}
\newtagform{normalsize}[\normalsize]{\normalsize(}{\normalsize)}
\usepackage{float}
\usepackage{natbib}
\usepackage[normalem]{ulem}
\usepackage{ifthen}
\usepackage[hidelinks]{hyperref}

\theoremstyle{plain}
\newtheorem{theorem}{Theorem}[section]
\newtheorem{lemma}[theorem]{Lemma}
\newtheorem{corollary}[theorem]{Corollary}
\newtheorem{proposition}[theorem]{Proposition}

\theoremstyle{definition}

\theoremstyle{remark}
\newtheorem{remark}[theorem]{Remark}

\providecommand{\keywords}[1]{\par\medskip\noindent\textbf{Keywords: }#1}

\begin{document}
\usetagform{normalsize} 

\title{Quantifying the Uncertainty of Uncertainty -  Extending Greenwood's Formula to Variances}

\author{J. Rodenkirchen\thanks{Corresponding author. Tel.: +49 521 106-86808. E-mail: \href{mailto:juergen.rodenkirchen@uni-bielefeld.de}{juergen.rodenkirchen@uni-bielefeld.de}.}\\[3pt]
\small Biostatistics and Medical Biometry, Medical School OWL,\\
\small Bielefeld University, Universit\"atsstr. 25, Bielefeld, 33615, Germany}
\hypersetup{pdfauthor={J. Rodenkirchen}}
\date{}
\maketitle

\begin{abstract}
Motivated by uncertainty quantification in survival analysis, we derive a
first-order estimator of the original-scale variance approximation for Greenwood's variance estimator at a fixed time point. The derivation uses the natural \(n\)-scaling of
Greenwood's estimator and of the risk-set sum entering Greenwood's formula and
is translated back to the original Greenwood scale. A directly computable
diagonal base contribution is obtained first, but this term alone does not
provide the full first-order variance. Two additional contributions are
required: a between-component covariance correction associated with the
Kaplan--Meier factor and the risk-set sum, and a pairwise correction for the remaining first-order variance of the Greenwood component. Both corrections are constructed from empirical cross-product blocks of event and at-risk indicators using the Delta method; the accompanying theory paper supplies the global first-order theory and the subsequent benchmark. Combining these terms yields an explicit count-based corrected estimator. Under independent right censoring and standard
regularity conditions in a continuous-time no-ties model, at a fixed time with positive limiting scaled Greenwood component $h(t)>0$, the $n^3$-scaled corrected estimator converges in probability to the limiting CLT variance from the joint first-order theory of the Kaplan--Meier estimator and the scaled Greenwood component. A separate moment-convergence argument identifies this limit with the limit of the $n^3$-scaled actual sampling variance, yielding consistency on that scale and relative consistency when the limiting variance is positive.
\end{abstract}

\section{Introduction}
\label{intro}

Precision medicine aims to identify treatment strategies that account for
heterogeneity in individual treatment response. Two central requirements are
therefore the presence of treatment-effect heterogeneity and the availability
of clinical predictors that identify patient groups with different expected
treatment benefits \citep{Den20,Chu20,Nol22,Kro22}.
Variability-based meta-regression approaches provide one possible way to
investigate such heterogeneity by comparing outcome variability between
treatment and control arms of randomized trials \citep{Kus23,Var24}.
Extending related ideas to time-to-event outcomes motivates the development of
non-parametric methods for quantifying uncertainty in survival analyses.

The original precision-medicine motivation concerns
variability in individual time-to-event outcomes, whereas Greenwood's formula
quantifies the sampling variability of the Kaplan--Meier estimator at a fixed
time point. These are distinct, although related, inferential problems. The
present article addresses the latter, namely the sampling variance of
Greenwood's variance estimator; extensions to variability measures defined
directly on the survival-time distribution are left for future work.

The Kaplan--Meier product-limit estimator \citep{Kap58} is the classical
non-parametric estimator of the survival function under right censoring, and
Greenwood's formula \citep{Gre26} provides its classical pointwise variance
estimator. The large-sample behavior of the Kaplan--Meier estimator and related
variance quantities has been studied extensively; classical contributions
include \citet{BreslowCrowley1974} and \citet{Gill1980,Gill1983}, with standard
counting-process treatments given, for example, by \citet{Fleming1991} and
\citet{Andersen1993}. Influence-function methods provide a complementary route
to asymptotic variance calculations for estimators constructed from the
Kaplan--Meier estimator; an early contribution in this direction is
\citet{Rei81}.

The specific first-order asymptotic theory required here is developed
separately in \citet{GHTheory26}. There, the appropriately scaled
risk-set sum entering Greenwood's formula is shown to admit an asymptotically
linear representation, and its joint first-order behavior with the
Kaplan--Meier estimator yields the asymptotic distribution of the scaled
Greenwood estimator. The present article addresses a complementary problem:
can the sampling variance of Greenwood's variance estimator itself be estimated
by an explicit formula based directly on the observed event counts and risk
sets?

The first-order theory developed in \citet{GHTheory26} identifies the limiting
CLT variance \(v_G(t)\) of the scaled, centered Greenwood error, determining the original-scale
first-order approximation \(v_G(t)/n^3\). This theoretical variance is not itself a
sample-based estimator. Rather than constructing an estimator by direct plug-in
into the full limiting variance formula, we pursue a constructive, count-based approach.
Using Greenwood's product representation and recursive finite-dimensional
Delta-method calculations, we first obtain a diagonal base contribution and
then construct the omitted between-component covariance correction and the
pairwise correction for the remaining first-order variance of the Greenwood component. The
resulting corrected estimator is an explicit formula in observed event counts
and risk-set sizes. Both corrections are constructed from empirical cross-product blocks; their count-based limits are obtained before comparison with the accompanying theory paper.

Under independent right censoring and standard regularity conditions in a continuous-time no-ties model, with $h(t)>0$ at the fixed time of interest, we show that the corrected estimator consistently estimates the original-scale first-order variance approximation on the $n^3$-scale, as validated against the limiting CLT variance from the joint first-order theory of the Kaplan--Meier estimator and the scaled Greenwood component. Thus, the diagonal
base formula serves as an intermediate decomposition, whereas the
covariance-corrected expression is the main inferential result. We also examine
the no-censoring case and the asymptotic non-negativity of the corrected
estimator. A separate moment-convergence argument establishes consistency for the actual sampling variance on the $n^3$-scale and relative consistency when $v_G(t)>0$.

The article is organized as follows. We first introduce the Kaplan--Meier
estimator, Greenwood's variance estimator, the risk-set sum entering Greenwood's
formula, and the scaling conventions used throughout the paper. We then derive
the preliminary diagonal base formula and identify the first-order contributions
that it does not contain. Next, the between-component covariance correction and
the additional variance contribution of the risk-set sum are derived and
combined with the base term to obtain the corrected estimator. The resulting
formula is subsequently validated against the joint first-order
influence-function theory, and its behavior in the no-censoring case and its
asymptotic non-negativity are considered. Consistency for the actual sampling variance is then established. The scope of the results and directions
for future work are discussed in Section~\ref{sec:discussion}. Auxiliary
first-order plug-in results are collected in the Appendix.

\newpage

\section{Setup, Notation, and Assumptions}
\label{definitions}

\subsection{Observed data and counting-process notation}

Fix a finite time horizon $\tau>0$. Let $T\ge0$ denote the event time with survival function
\[
S(t)=P(T>t),
\]
and let $C\ge0$ denote the censoring time. We allow $C=\infty$ to represent no censoring. For a generic subject, define
\[
X=T\wedge C,
\qquad
\Delta=\mathbf 1\{T\le C\},
\qquad
Z=(X,\Delta).
\]
For subjects $i=1,\ldots,n$, define analogously
\[
X_i=T_i\wedge C_i,
\qquad
\Delta_i=\mathbf 1\{T_i\le C_i\},
\qquad
Z_i=(X_i,\Delta_i).
\]

For $t\le\tau$, let
\[
Y_i(t)=\mathbf 1\{X_i\ge t\},
\qquad
Y_n(t)=\sum_{i=1}^n Y_i(t)
\]
denote the individual and aggregate at-risk processes, and let
\[
N_i(t)=\mathbf 1\{X_i\le t,\Delta_i=1\},
\qquad
N_n(t)=\sum_{i=1}^n N_i(t)
\]
denote the corresponding failure counting processes. The population at-risk probability is
\[
y(t):=P(X\ge t)=E\{Y_i(t)\}.
\]

Let
\[
0<t_1<t_2<\cdots
\]
denote the distinct observed failure times. At $t_j$, write
\[
d_j:=\Delta N_n(t_j),
\qquad
n_j:=Y_n(t_j-),
\]
so that $d_j$ is the number of observed failures and $n_j$ is the number of individuals at risk immediately before $t_j$.
Let $D_i^C(u):=\mathbf{1}_{\{X_i\le u,\ \Delta_i=0\}}$ denote the observed censoring counting process. We use the usual completed, right-continuous augmentation of the filtration generated by the individual failure and censoring counting processes up to each time.
The at-risk processes are left-continuous and predictable with respect to the observed-history filtration; in particular, $Y_n(t_j-)=Y_n(t_j)$.

Under the assumptions stated below, these counting processes admit the standard multiplicative intensity and martingale representations used in the auxiliary asymptotic arguments. The corresponding martingale notation is introduced only where it is needed.

\subsection{Kaplan--Meier and Greenwood quantities}

The Kaplan--Meier product-limit estimator \citep{Kap58} is
\begin{equation}
\label{eq:KaplanMeier}
\widehat S_n(t)
=
\prod_{j:t_j\le t}
\left(1-\frac{d_j}{n_j}\right).
\end{equation}

Greenwood's formula \citep{Gre26} estimates the pointwise sampling variance of $\widehat S_n(t)$. For convenience, we refer to the risk-set sum
\begin{equation}
\label{eq:Hhat}
\widehat H_n(t)
=
\sum_{j:t_j\le t}
\frac{d_j}{n_j(n_j-d_j)}
\end{equation}
as the \emph{Greenwood component}. Greenwood's variance estimator can then be written as
\begin{equation}
\label{eq:Greenwood}
\widehat G_n(t)
=
\widehat S_n(t)^2\widehat H_n(t).
\end{equation}

\subsection{Global well-definedness}
\label{subsec:global_well_definedness}

Let $n\ge1$. The usual Greenwood expression contains denominators $n_j-d_j$, which need not be strictly positive for every finite sample. To make the Greenwood quantities globally defined as random variables, let
\begin{equation}
\label{eq:Dn_def}
\mathcal D_n(t)
:=
\left\{
n_j-d_j>0
\ \text{for every observed failure time }t_j\le t
\right\}.
\end{equation}
For $n\ge2$, on $\mathcal D_n(t)$, $\widehat H_n(t)$ is given by \eqref{eq:Hhat}. For $n<2$ or on $\mathcal D_n(t)^c$, set $\widehat H_n(t)=0$ and set all later Greenwood-denominator-dependent count quantities and variance contributions to zero, including $\widehat Q_n$, $\widehat C_j$, $\widehat C_{\mathrm{diag}}$, $\widehat B_{\mathrm{base}}$, the empirical correction terms, and $\widehat R_{\mathrm{base}}$, $\widehat R_{\mathrm{corr}}$. The Kaplan--Meier estimator retains its usual product definition, and $\widehat G_n(t)=\widehat S_n(t)^2\widehat H_n(t)$ everywhere.

Throughout, the displayed count-based formulas are interpreted for $n\ge2$ on $\mathcal D_n(t)$ unless otherwise specified. On this regular set with no observed failures, empty sums are zero, $\widehat S_n(t)=1$, and all empirical quantities involving division by $\widehat H_n(t)$ are set to zero. In particular, the definitions by quotient and by expanded sums agree under these conventions. Pointwise auxiliary weights in the Appendix are extended separately.

Under the assumptions below,
\begin{equation}
\label{eq:Dn_prob}
P\{\mathcal D_n(t)\}\longrightarrow1
\qquad
\text{for every fixed }t\le\tau.
\end{equation}
Indeed, under the no-ties model $d_j=1$, while
\[
\frac{Y_n(\tau)}{n}
\xrightarrow{\mathrm{a.s.}}
y(\tau)>0.
\]
Hence $P\{Y_n(\tau)\ge2\}\to1$, and on this event
\[
n_j\ge Y_n(\tau)\ge2
\]
for every observed failure time $t_j\le\tau$. Thus the fixed zero extension on $\mathcal D_n(t)^c$ is asymptotically irrelevant.

This denominator condition is separate from positivity required by logarithmic transformations. Whenever log-scale Delta-method calculations are used below, we additionally restrict attention to fixed time points satisfying
\begin{equation}
\label{eq:h_positive_log}
h(t)>0,
\end{equation}
where $h(t)$ is defined below. Consistency then implies that the corresponding Greenwood quantities are positive with probability tending to one, so that their logarithms are asymptotically well defined. They may again be extended arbitrarily on the asymptotically negligible complementary event.

\subsection{Scaling and pointwise convention}

The quantities satisfy $\widehat H_n(t)=O_p(n^{-1})$ and $\widehat G_n(t)=O_p(n^{-1})$. Their natural $O_p(1)$-scaled versions are
\begin{equation}
\label{eq:scaled_greenwood_quantities}
\widetilde H_n(t)
:=
n\widehat H_n(t),
\qquad
\widetilde G_n(t)
:=
n\widehat G_n(t)
=
\widehat S_n(t)^2\widetilde H_n(t).
\end{equation}

Throughout this manuscript, asymptotic statements are pointwise in time unless stated otherwise; that is, $t\le\tau$ is fixed as $n\to\infty$.

\subsection{Original-scale first-order variance and covariance convention}
\label{subsec:Avar_convention}

Because several quantities in this paper live on different asymptotic scales,
we distinguish variances and covariances appearing in asymptotic limits from
the corresponding first-order approximations on the original scale.

Suppose that, for deterministic centering constants $\theta_n$ and scaling
constants $a_n\to\infty$,
\[
a_n\{X_n-\theta_n\}
\xRightarrow{d}
N(0,\sigma_X^2).
\]
We refer to $\sigma_X^2$ as the \emph{asymptotic variance of the scaled
error}. The corresponding first-order variance approximation for $X_n$ on
its original scale is
\[
A.Var(X_n)
:=
\frac{\sigma_X^2}{a_n^2}.
\]

Analogously, if
\[
a_n
\begin{pmatrix}
X_n-\theta_{X,n}\\
Y_n-\theta_{Y,n}
\end{pmatrix}
\xRightarrow{d}
N_2\left\{
\begin{pmatrix}
0\\
0
\end{pmatrix},
\begin{pmatrix}
\sigma_X^2 & \sigma_{XY}\\
\sigma_{XY} & \sigma_Y^2
\end{pmatrix}
\right\},
\]
we define the corresponding original-scale first-order covariance by
\[
A.Cov(X_n,Y_n)
:=
\frac{\sigma_{XY}}{a_n^2}.
\]
The same convention is used componentwise for first-order covariance
matrices.

Thus the limiting variances and covariances describe the joint limiting
distribution after scaling, whereas $A.Var$ and $A.Cov$ translate these
quantities back to the scale on which the estimators themselves are
reported. They are first-order asymptotic approximations and do not by
themselves assert convergence of the corresponding exact finite-sample
variances or covariances. In particular, the first-order validation alone does not establish $n^3\operatorname{Var}\{\widehat G_n(t)\}\to v_G(t)$; this follows from the separate moment-convergence result in Proposition~\ref{prop:G_moment_convergence}.

For an ordinary $\sqrt n$-estimator,
\[
A.Var(X_n)
=
\frac{\sigma_X^2}{n},
\]
and analogously $A.Cov(X_n,Y_n)=\sigma_{XY}/n$ under a joint
$\sqrt n$ limit.

Greenwood's estimator $\widehat G_n(t)$, however, is $O_p(n^{-1})$. Its regular first-order analysis is therefore carried out for the
scaled quantity
\[
\widetilde G_n(t)=n\widehat G_n(t),
\]
and the resulting asymptotic variance is subsequently translated back to
the original scale of $\widehat G_n(t)$.

\subsection{Regularity assumptions}

The asymptotic results are developed under the following conditions.

\begin{enumerate}[label=(A\arabic*), leftmargin=2.4em]

\item
\label{ass:iid-independent-censoring}
\textbf{Sampling and independent censoring.}
The pairs $(T_i,C_i)$, $i=1,\ldots,n$, are i.i.d., and
\[
T_i\perp C_i.
\]

\item
\label{ass:positivity}
\textbf{Positive at-risk probability.}
With $y(t)=P(X\ge t)$,
\[
\inf_{t\in[0,\tau]}y(t)>0.
\]

\item
\label{ass:hazard-regularity}
\textbf{Event-hazard regularity.}
The event-time distribution is absolutely continuous on $[0,\tau]$, with hazard function $\lambda$, and
\[
\Lambda(t)
:=
\int_0^t\lambda(u)\,du
<
\infty,
\qquad
t\le\tau.
\]
No continuity or hazard-function assumption on the censoring-time distribution is required beyond the independent-censoring condition in (A1).

\item
\label{ass:ties}
\textbf{Continuous-time no-ties model.}
Observed failure times are almost surely distinct, so that
\[
d_j=1
\]
at every observed failure time. The general $d_j$-notation is nevertheless retained where algebraically convenient because it represents Greenwood's exact finite-sample formula naturally. The asymptotic validation of the full corrected estimator in this paper, however, is carried out under the continuous-time no-ties model.

\end{enumerate}

\subsection{First-order theory used in this paper}

Consistency of the scaled Greenwood component is already contained in
\citet[Theorem~4.2.2]{Gill1980}. A direct proof under the present assumptions
is also given in \citet[Theorem~5]{GHTheory26}. The latter paper develops the
first-order fluctuation theory used below. Under
Assumptions~\ref{ass:iid-independent-censoring}--\ref{ass:ties},
for every fixed $t\le\tau$,
\[
\widetilde H_n(t)
=
n\widehat H_n(t)
\xrightarrow{P}
h(t),
\qquad
h(t)
:=
\int_0^t
\frac{\lambda(u)}{y(u)}
\,du,
\]
and
\[
\sqrt n
\left\{
\widetilde H_n(t)-h(t)
\right\}
\xRightarrow{d}
N\{0,v_H(t)\}.
\]
The second assertion follows from Theorem~7 of \citet{GHTheory26}.

Together with the standard first-order theory of the Kaplan--Meier estimator,
Theorem~9 of \citet{GHTheory26} further yields the joint central limit theorem
\[
\sqrt n
\begin{pmatrix}
\widehat S_n(t)-S(t)\\
\widetilde H_n(t)-h(t)
\end{pmatrix}
\xRightarrow{d}
N_2\{0,\Sigma(t)\},
\]
so that the multivariate Delta method applies to
\[
\widetilde G_n(t)
=
\widehat S_n(t)^2\widetilde H_n(t).
\]

The separate theory paper thus supplies the global first-order theory and a benchmark for the original-scale variance approximation of Greenwood's variance estimator. The present paper addresses the complementary problem of constructing an explicit count-based estimator on the original Greenwood scale whose scaled limit reproduces that asymptotic variance.

\section{Base Variance Formula for Greenwood's Estimator}
\label{sec:main}

Consider Greenwood's estimator
\[
\widehat G_n(t)=\widehat S_n(t)^2\widehat H_n(t)
\]
from \eqref{eq:Greenwood}. The original-scale first-order variance approximation of the unscaled Greenwood estimator $\widehat G_n(t)$ contains several distinct contributions. We first isolate a directly computable diagonal base
contribution on the original Greenwood scale. This base contribution is an
intermediate component of the derivation and is not, by itself, the full
first-order variance estimator.

We denote the resulting plug-in base estimator by
\(\widehat R_{\mathrm{base}}(t)\). It retains the contribution associated with
the Kaplan--Meier factor together with the diagonal martingale contribution
to the variance of the Greenwood component. It does not yet include the
between-component covariance correction \(\Gamma_{\widehat A}(t)\) or the
additional Greenwood-component variance contribution
\(\Gamma_{\widehat B}(t)\). These terms are incorporated in
Section~\ref{sec:improvement_base_formula}, where the corrected estimator
\(\widehat R_{\mathrm{corr}}(t)\) is obtained. The corrected estimator, rather
than the base formula stated here, is the main inferential result of the
paper.

The base result is as follows.

\begin{theorem}[Base variance formula for Greenwood's estimator]
\label{thm:main}
Assume Assumptions~\ref{ass:iid-independent-censoring}--%
\ref{ass:ties} and the log-positivity condition~\eqref{eq:h_positive_log} for the fixed time point $t\le \tau$.
Let $\widehat S_n$ be the Kaplan--Meier estimator \eqref{eq:KaplanMeier},
and let
\[
\widehat G_n(t)=\widehat S_n(t)^2\widehat H_n(t),
\qquad
\widehat H_n(t)
=
\sum_{j:t_j\le t}
\frac{d_j}{n_j(n_j-d_j)}
\]
be Greenwood's estimator and its unscaled Greenwood component. Then the
directly computable plug-in estimator of the diagonal base contribution to
the original-scale first-order variance approximation of $\widehat G_n(t)$ is
\begin{equation}
\label{eq:hatR}
\boxed{
\widehat R_{\mathrm{base}}(t)
=
\widehat S_n(t)^4
\left\{
4\widehat H_n(t)^3
+
\sum_{j:t_j\le t}
\frac{d_j}{n_j(n_j-d_j)^3}
\right\}.
}
\end{equation}
Moreover, define
\[
L(t)
:=
\int_0^t\frac{\lambda(u)}{y(u)^3}\,du.
\]
Then the base estimator satisfies
\[
n^3\widehat R_{\mathrm{base}}(t)
\xrightarrow{P}
S(t)^4\{4h(t)^3+L(t)\}.
\]
Thus the base formula has a precise first-order asymptotic target, but this
target is not in general the full limiting CLT variance of the scaled, centered Greenwood error. The missing first-order contributions are incorporated through
\(\Gamma_{\widehat A}(t)\) and \(\Gamma_{\widehat B}(t)\) in
Section~\ref{sec:improvement_base_formula}.

\end{theorem}

\section{First-Order Delta-Method Framework}
\label{Delta-Method}

We use first-order Delta-method relations interpreted according to the
original-scale variance and covariance convention fixed in
Section~\ref{subsec:Avar_convention}. The relations below concern the
theoretical first-order quantities $A.Var$ and $A.Cov$. When they are used
to construct an estimator, unknown population quantities are subsequently
replaced by their corresponding sample-based plug-in quantities. For a
detailed discussion of the Delta method, see, for example, \citep{Rao73}.

\subsection{Univariate first-order variance propagation}

Let \(X_n\) be a regular estimator of a scalar parameter \(\theta\), and
suppose that
\[
\sqrt n\,(X_n-\theta)
\xRightarrow{d}
N(0,\sigma_X^2).
\]
Then the limiting CLT variance is \(\sigma_X^2\), whereas, according to
Section~\ref{subsec:Avar_convention}, the corresponding original-scale
first-order variance approximation is
\[
A.Var(X_n)=\frac{\sigma_X^2}{n}.
\]
If \(g\) is differentiable at \(\theta\), then the Delta method gives
\[
\sqrt n\,\{g(X_n)-g(\theta)\}
\xRightarrow{d}
N\left(0,\{g'(\theta)\}^2\sigma_X^2\right),
\]
and hence, on the original scale,
\begin{equation}
\label{eq:delta_var_original}
A.Var\{g(X_n)\}
=
\{g'(\theta)\}^2 A.Var(X_n).
\end{equation}
Equation~\eqref{eq:delta_var_original} is a theoretical first-order
relation. When it is used to construct a sample-based variance estimator,
unknown quantities such as \(\theta\) and $A.Var(X_n)$ are replaced by
appropriate plug-in estimators.

\subsection{First-order covariance propagation}

If \(X_n\) and \(Y_n\) are jointly regular estimators of scalar parameters
\(\theta_X\) and \(\theta_Y\), and if \(f\) and \(g\) are differentiable at
these parameters, then the corresponding original-scale first-order
covariance relation is
\begin{equation}
\label{eq:delta_cov_original}
A.Cov\{f(X_n),g(Y_n)\}
=
f'(\theta_X)g'(\theta_Y)\,A.Cov(X_n,Y_n).
\end{equation}
Equivalently, for vector-valued jointly regular estimators, the covariance
propagation takes the form
\begin{equation}
\label{eq:delta_cov_multivariate}
A.Cov\{f(\mathbf X_n),g(\mathbf Y_n)\}
=
\nabla f(\boldsymbol\theta_X)^\top
A.Cov(\mathbf X_n,\mathbf Y_n)
\nabla g(\boldsymbol\theta_Y).
\end{equation}

As in the variance case, these are theoretical first-order relations;
sample-based formulas arise only after the unknown population quantities
are replaced by appropriate plug-in estimators. In the Greenwood
derivation below, the between-component contribution
\(\Gamma_{\widehat A}(t)\) enters through such first-order covariance
propagation. The additional Greenwood-component variance contribution
\(\Gamma_{\widehat B}(t)\) is distinct and completes the first-order
variance of the Greenwood component beyond the diagonal base contribution.
Both terms are incorporated in Section~\ref{sec:improvement_base_formula}.

\section{Derivation of the Base Variance Formula}
\label{proof}

We now derive the count-based base estimator in
Theorem~\ref{thm:main}. The derivation proceeds on the log scale, where the
product structure of Greenwood's estimator becomes additive. We first
establish the theoretical first-order back-transformation from the log
scale to the original Greenwood scale. We then decompose the log-scale
first-order variance and construct plug-in estimators for the contributions
retained in the base formula. The between-component covariance contribution
\(\Gamma_{\widehat A}(t)\) and the additional Greenwood-component variance
contribution \(\Gamma_{\widehat B}(t)\) are incorporated subsequently in
Section~\ref{sec:improvement_base_formula}.

\subsection{Log-Scale Back-Transformation}
\label{4.1}

The regular first-order object associated with Greenwood's estimator is
\[
\widetilde G_n(t)=n\widehat G_n(t).
\]
Corollary~11 of \citet{GHTheory26} gives
\[
\sqrt n\{
\widetilde G_n(t)-S(t)^2h(t)
\}
\xRightarrow{d}
N\{0,v_G(t)\},
\]
where \(v_G(t)\) denotes the variance of the limiting normal distribution. Hence, according to the original-scale convention of
Section~\ref{subsec:Avar_convention},
\[
A.Var\{\widehat G_n(t)\}
=
\frac{v_G(t)}{n^3}.
\]
Here, \(A.Var\) denotes the original-scale first-order variance approximation as defined in
Section~\ref{subsec:Avar_convention}, rather than the exact finite-sample variance.

Because \(h(t)>0\), the log transformation is asymptotically well defined.
Moreover,
\[
\log\widehat G_n(t)
=
\log\widetilde G_n(t)-\log n.
\]
Applying the Delta method to the regular scaled quantity
\(\widetilde G_n(t)\) therefore gives
\[
\sqrt n
\left[
\log\widehat G_n(t)
-
\log\left\{\frac{S(t)^2h(t)}{n}\right\}
\right]
\xRightarrow{d}
N\left\{
0,
\frac{v_G(t)}{S(t)^4h(t)^2}
\right\}.
\]
Consequently,
\[
A.Var\{\log\widehat G_n(t)\}
=
\frac{v_G(t)}
     {nS(t)^4h(t)^2}.
\]
Combining the preceding two relations yields the theoretical first-order
identity
\begin{equation}
\label{eq:base}
A.Var\{\widehat G_n(t)\}
=
\left\{\frac{S(t)^2h(t)}{n}\right\}^2
A.Var\{\log\widehat G_n(t)\}.
\end{equation}

The factor \(S(t)^2h(t)/n\) in \eqref{eq:base} is consistently estimated by
\[
\widehat G_n(t)
=
\widehat S_n(t)^2\widehat H_n(t).
\]
More precisely, since $h(t)>0$,
\[
\frac{\widehat G_n(t)}{S(t)^2h(t)/n}\xrightarrow{P}1.
\]
Thus, when a sample-based estimator of the log-scale first-order variance
has been constructed, the corresponding plug-in back-transformation to the
original Greenwood scale uses the factor \(\widehat G_n(t)^2\).

\subsection{Log-Scale Variance Decomposition}
\label{4.2}

From
\[
\log\widehat G_n(t)
=
2\log\widehat S_n(t)
+
\log\widehat H_n(t),
\]
the first-order variance decomposition is
\[
\begin{aligned}
A.Var\{\log\widehat G_n(t)\}
={}&
4A.Var\{\log\widehat S_n(t)\}
+
A.Var\{\log\widehat H_n(t)\} \\
&+
\Gamma_{\widehat A}(t),
\end{aligned}
\]
where
\begin{equation}
\label{eq:Gamma_hat_A_1}
\Gamma_{\widehat A}(t)
=
4A.Cov\{
\log\widehat S_n(t),
\log\widehat H_n(t)
\}.
\end{equation}

The first-order log-scale variance of the Kaplan--Meier estimator (Appendix Lemma~\ref{lem:greenwood_log_variance_base}) is
\[
A.Var\{\log\widehat S_n(t)\}
=
\frac{h(t)}{n}.
\]
The required consistency is already contained in
\citet[Theorem~4.2.2]{Gill1980} and is also proved directly in
\citet[Theorem~5]{GHTheory26}:
\[
n\widehat H_n(t)
\xrightarrow{P}
h(t).
\]
Thus the Greenwood component \(\widehat H_n(t)\) is the corresponding
sample-based plug-in estimator.

The base formula retains this Kaplan--Meier-factor contribution but leaves
the between-component term \(\Gamma_{\widehat A}(t)\) outside the base
approximation. It also retains only the diagonal base part of
\(A.Var\{\log\widehat H_n(t)\}\), which is constructed next.

\subsection{Diagonal Base Contribution and Deferred Pairwise Terms}
\label{4.3}

The remaining log-scale first-order term from Subsection~\ref{4.2} is
\(A.Var\{\log\widehat H_n(t)\}\). Theorem~7 of
\citet{GHTheory26} gives the central limit theorem for
\(\widetilde H_n(t)=n\widehat H_n(t)\):
\[
\sqrt n\{\widetilde H_n(t)-h(t)\}
\xRightarrow{d} N\{0,v_H(t)\}.
\]
For \(h(t)>0\), applying the log-Delta method to the scaled quantity $\widetilde H_n(t)$ and using $\log\widehat H_n(t)=\log\widetilde H_n(t)-\log n$ on the positive set gives
\[
\begin{aligned}
A.Var\{\log\widehat H_n(t)\}
&= \frac{1}{\{h(t)/n\}^2}\,A.Var\{\widehat H_n(t)\} \\
&= \frac{1}{\{h(t)/n\}^2}\,\frac{v_H(t)}{n^3}
 = \frac{v_H(t)}{n\,h(t)^2}.
\end{aligned}
\]
We now construct contributions to the original-scale variance of
\(\widehat H_n(t)\), before applying the displayed log-Delta factor
through its sample-based plug-in.

The unscaled Greenwood component is the sum of its event-time contributions,
\[
\widehat H_n(t)
=
\sum_{j:t_j\le t}D_j,
\qquad
D_j
:=
\frac{d_j}{n_j(n_j-d_j)}.
\]
For a fixed index set, the familiar variance identity reads
\[
\operatorname{Var}\left(\sum_{j=1}^{m}D_j\right)
=
\sum_{j=1}^{m}\operatorname{Var}(D_j)
+
\Gamma_{\widehat B}^{(m)},
\]
where the off-diagonal Greenwood-summand covariance contribution is
explicitly
\begin{equation}
\label{eq:GammaB_structural}
\Gamma_{\widehat B}^{(m)}
:=2\sum_{1\le i<j\le m}\operatorname{Cov}(D_i,D_j).
\end{equation}
For the observed event-time sum up to \(t\), we write
\(\Gamma_{\widehat B}(t)\) for the corresponding \emph{first-order
off-diagonal contribution}. Since the event times and the number of
summands are random, the fixed-index identity is not asserted as an
exact unconditional variance identity for \(\widehat H_n(t)\).
Instead, it motivates the count-based pairwise Delta-method construction in
Subsection~\ref{subsec:GammaB_construction}.
The term off-diagonal refers to the pairwise construction of the empirical candidate. Its $n^3$-scaled first-order target is the remainder of the Greenwood-component CLT variance after subtracting the diagonal martingale contribution, as identified in Appendix Lemma~\ref{lem:GammaB_limit}. No covariance-sum representation for the actual randomly indexed summands is asserted.

The diagonal contribution is derived in Subsection~\ref{4.4}; the
first-order limits of the diagonal and pairwise constructions are
established in Appendix Lemmas~\ref{lem:diag_limit}
and~\ref{lem:GammaB_limit}, respectively.

The local event-count Delta calculation detailed in
Subsection~\ref{4.4} yields the directly computable diagonal quantity
\begin{equation}
\label{eq:Cdiag_def}
\widehat C_{\mathrm{diag}}(t)
:=
\sum_{j:t_j\le t}
\frac{d_j}{n_j(n_j-d_j)^3}.
\end{equation}
Under the continuous-time no-ties asymptotic model,
\[
\widehat C_{\mathrm{diag}}(t)
=
\sum_{j:t_j\le t}
\frac{1}{n_j(n_j-1)^3},
\]
and Appendix Lemma~\ref{lem:diag_limit} gives
\begin{equation}
\label{eq:Cdiag_limit_main}
n^3\widehat C_{\mathrm{diag}}(t)
\xrightarrow{P}
L(t).
\end{equation}

Using the log-Delta relation established at the start of this subsection,
we transfer the count-based diagonal quantity to the log scale.
Accordingly, define
\begin{equation}
\label{eq:hatB3}
\widehat B_{\mathrm{base}}(t)
:=
\widehat H_n(t)^{-2}
\widehat C_{\mathrm{diag}}(t)
\qquad (\widehat H_n(t)>0).
\end{equation}
Set $\widehat B_{\mathrm{base}}(t)=0$ whenever $\widehat H_n(t)=0$, using the global extension convention of Subsection~\ref{subsec:global_well_definedness} on $\mathcal D_n(t)^c$.
Since \(n\widehat H_n(t)\xrightarrow{P}h(t)\) with $h(t)>0$, the following quotient calculation is valid with probability tending to one:
\[
n\widehat B_{\mathrm{base}}(t)
=
\frac{n^3\widehat C_{\mathrm{diag}}(t)}
     {\{n\widehat H_n(t)\}^2}
\xrightarrow{P}
\frac{L(t)}{h(t)^2}.
\]
Thus \(\widehat B_{\mathrm{base}}(t)\) estimates \(L(t)/\{n h(t)^2\}\)
to first order.
The limit validates the aggregated count-based plug-in at first order; it is not used to construct its count-based formula.

\subsection{Local Delta Representation of the Diagonal Plug-In}
\label{4.4}

We now verify the local Delta-method calculation underlying the
count-based diagonal plug-in introduced in \eqref{eq:Cdiag_def}. Write
\[
\widehat D_j
=
\frac{d_j}{n_j(n_j-d_j)}
=
\frac{\widehat\lambda_j}{n_j(1-\widehat\lambda_j)},
\qquad
\widehat\lambda_j=\frac{d_j}{n_j}.
\]
Here $\widehat\lambda_j$ is the observed event fraction, not an instantaneous hazard rate.
Holding the predictable risk-set size \(n_j\) fixed for this local
calculation, let
\[
f_j(x)=\frac{x}{n_j(1-x)},
\qquad
f_j'(x)=\frac{1}{n_j(1-x)^2}.
\]
The local event-count Delta calculation gives
\[
\left\{f_j'(\widehat\lambda_j)\right\}^2
\frac{\widehat\lambda_j(1-\widehat\lambda_j)}{n_j}
=
\frac{d_j}{n_j(n_j-d_j)^3}.
\]
Thus define
\begin{equation}
\label{eq:hatCj2}
\widehat C_j(t)
:=
\frac{d_j}{n_j(n_j-d_j)^3},
\end{equation}
so that
\[
\widehat C_{\mathrm{diag}}(t)
=
\sum_{j:t_j\le t}\widehat C_j(t).
\]

This local calculation is the Delta-method construction of the
retained diagonal term. Its asymptotic justification in the present
paper is \eqref{eq:Cdiag_limit_main}, not an assumption that the random
Greenwood summands are independent or that the risk sets are nonrandom.

\begin{remark}[Interpretation of the local Delta-method construction]
Under the continuous-time no-ties assumption, the counting-process jumps satisfy $\Delta N_n(u)\in\{0,1\}$. Since $t_j$ indexes an observed failure time, $d_j=\Delta N_n(t_j)=1$ whenever this time exists. Holding the risk-set size $n_j$ fixed therefore leaves no randomness in the observed event fraction $d_j/n_j$ at that existing failure time. Consequently, the binomial-form factor used above is not interpreted as an actual conditional variance at an observed failure time, but as a local plug-in motivating the count-based Delta-method construction. Its first-order validity is established for the aggregated diagonal contribution through \eqref{eq:Cdiag_limit_main}.
\end{remark}

\begin{remark}[Relation to martingale variation]
\nopagebreak[4]
Under the no-ties assumptions, let $w_G(m)=1/[m(m-1)]$ for $m\ge2$ and $w_G(m)=0$ otherwise, and set $M_H(t)=\int_0^t w_G(Y_n(u))\,dM_n(u)$, where $M_n(u)=N_n(u)-\int_0^u Y_n(s)\lambda(s)\,ds$ is the counting-process martingale. The integrand is bounded and predictable, and the compensator of $N_n$ is bounded by $n\Lambda(t)$ on $[0,t]$; hence $M_H$ is square-integrable. For $n\ge2$ on $\mathcal D_n(t)$,
\[
[M_H](t)=\sum_{j:t_j\le t}\frac{1}{n_j^2(n_j-1)^2},\qquad
\widehat C_{\mathrm{diag}}(t)-[M_H](t)
=\sum_{j:t_j\le t}\frac{1}{n_j^2(n_j-1)^3}.
\]
Uniform risk-set positivity and at most $n$ failures give a difference of order $O_p(n^{-4})$. Furthermore,
\[
n^3\langle M_H\rangle(t)
=n^3\int_0^t w_G(Y_n(u))^2Y_n(u)\lambda(u)\,du
\xrightarrow{P}L(t),
\]
by uniform risk-set convergence on the high-probability event $\inf_{u\le t}Y_n(u)/n\ge c$ for any fixed $0<c<y(t)$. Thus the diagonal plug-in shares the first-order limit of the optional and predictable quadratic variations of this martingale. This observation does not identify it with the actual sampling variance of the Greenwood component or by itself establish moment convergence for Greenwood's estimator; the latter is proved separately in Proposition~\ref{prop:G_moment_convergence}.
\end{remark}

\begin{remark}[Asymptotic interpretation of the deferred correction]
\label{rem:asymptotic_orientation_B}
The diagonal plug-in constructed above satisfies
\[
n^3\widehat C_{\mathrm{diag}}(t)\xrightarrow{P}L(t)
\]
by Appendix Lemma~\ref{lem:diag_limit}. The deferred pairwise correction
\(\widehat\Gamma_{\widehat B,\mathrm{nt}}(t)\), constructed by the bivariate
Delta method in Subsection~\ref{subsec:GammaB_construction}, has the limit
\[
n^3\widehat\Gamma_{\widehat B,\mathrm{nt}}(t)\xrightarrow{P}K_B(t),
\]
where $K_B(t)$ is the integral in \eqref{eq:GammaB_K_integral}, derived directly from the empirical pair kernel in Appendix Lemma~\ref{lem:GammaB_limit}. The two limits are compared with the full limiting variance $v_H(t)$ in Subsection~\ref{subsec:Rcorr_validation}.
\end{remark}

\subsection{Aggregation of the Base Contributions}
\label{4.5}

The sample-based base contribution to the log-scale variance of Greenwood's
estimator is
\[
\begin{aligned}
\widehat A_{\mathrm{base}}(t)
&:=4\widehat H_n(t)+\widehat B_{\mathrm{base}}(t)\\
&=4\widehat H_n(t)+\widehat H_n(t)^{-2}\widehat C_{\mathrm{diag}}(t),\qquad \widehat H_n(t)>0.
\end{aligned}
\]
The displayed quotient is used only on the positive set; the first line defines $\widehat A_{\mathrm{base}}(t)$ globally. The two consistency relations established above imply
\[
n\widehat A_{\mathrm{base}}(t)
\xrightarrow{P}
4h(t)+\frac{L(t)}{h(t)^2}.
\]

Using the plug-in back-transformation motivated by
\eqref{eq:base}, define
\[
\widehat R_{\mathrm{base}}(t)
:=
\widehat G_n(t)^2
\widehat A_{\mathrm{base}}(t).
\]
Since
\[
\widehat G_n(t)
=
\widehat S_n(t)^2\widehat H_n(t),
\]
we obtain, on the regular set with $\widehat H_n(t)>0$,
\[
\begin{aligned}
\widehat R_{\mathrm{base}}(t)
&=
\widehat S_n(t)^4
\widehat H_n(t)^2
\left\{
4\widehat H_n(t)
+
\widehat H_n(t)^{-2}
\widehat C_{\mathrm{diag}}(t)
\right\} \\
&=
\widehat S_n(t)^4
\left\{
4\widehat H_n(t)^3
+
\widehat C_{\mathrm{diag}}(t)
\right\} \\
&=
\widehat S_n(t)^4
\left\{
4\widehat H_n(t)^3
+
\sum_{j:t_j\le t}
\frac{d_j}{n_j(n_j-d_j)^3}
\right\},
\end{aligned}
\]
which is exactly \eqref{eq:hatR}. Under the common zero-extension convention, the expanded identity also holds on the remaining samples.

Finally,
\[
n\widehat G_n(t)
\xrightarrow{P}
S(t)^2h(t)
\]
and
\[
n\widehat A_{\mathrm{base}}(t)
\xrightarrow{P}
4h(t)+\frac{L(t)}{h(t)^2}.
\]
Therefore
\[
\begin{aligned}
n^3\widehat R_{\mathrm{base}}(t)
&=
\{n\widehat G_n(t)\}^2
\{n\widehat A_{\mathrm{base}}(t)\} \\
&\xrightarrow{P}
S(t)^4h(t)^2
\left\{
4h(t)+\frac{L(t)}{h(t)^2}
\right\} \\
&=
S(t)^4\{4h(t)^3+L(t)\}.
\end{aligned}
\]
This proves both the count-based base formula and its first-order
asymptotic characterization stated in Theorem~\ref{thm:main}.

\section{Full Corrected Variance Formula}
\label{sec:improvement_base_formula}

The base formula in Theorem~\ref{thm:main} combines the retained
Kaplan--Meier-factor contribution with the diagonal contribution to the
Greenwood-component variance. It does not yet contain the full first-order
between-component covariance contribution or the remaining first-order
Greenwood-component variance contribution. We now add these two terms and
obtain the corrected original-scale variance estimator.

For subject $r=1,\ldots,n$, let
\[
\begin{aligned}
E_{rj}&:=\mathbf1\{\text{subject }r\text{ has an observed failure at }t_j\},\\
R_{rj}&:=\mathbf1\{\text{subject }r\text{ is at risk just before }t_j\}.
\end{aligned}
\]
Write $W_{rj}=(E_{rj},R_{rj})^\top$ and $w_j=\sum_r W_{rj}=(d_j,n_j)^\top$. For $n\ge2$, define the empirical cross-product block
\[
C_{ij}:=\frac{n}{n-1}\sum_{r=1}^n
(W_{ri}-w_i/n)(W_{rj}-w_j/n)^\top
=\frac{n\sum_r W_{ri}W_{rj}^\top-w_iw_j^\top}{n-1}.
\]
Under the single-failure, nested-risk-set structure, for $i<j$,
\[
\begin{aligned}
\sum_r E_{ri}E_{rj}&=\sum_r E_{ri}R_{rj}=0,\\
\sum_r R_{ri}E_{rj}&=d_j,\qquad \sum_r R_{ri}R_{rj}=n_j.
\end{aligned}
\]
Consequently,
\[
C_{ij}=\frac{1}{n-1}
\begin{pmatrix}
-d_id_j & -d_i n_j\\
d_j(n-n_i) & n_j(n-n_i)
\end{pmatrix},\qquad i<j.
\]
For the diagonal blocks, $E_{ri}^2=E_{ri}$, $R_{ri}^2=R_{ri}$ and $E_{ri}R_{ri}=E_{ri}$, so
\[
C_{ii}=\frac{1}{n-1}
\begin{pmatrix}
d_i(n-d_i)&d_i(n-n_i)\\
d_i(n-n_i)&n_i(n-n_i)
\end{pmatrix},\qquad C_{ji}=C_{ij}^\top.
\]

\begin{remark}[Interpretation of the empirical cross-product blocks]
The no-ties assumption implies $\Delta N_n(u)\in\{0,1\}$, while the indexing by observed failure times implies $d_j=\Delta N_n(t_j)=1$ whenever the corresponding failure exists. The failure times, risk-set sizes and number of failures remain random. For $i<j$, conditional on $N_n(t)\ge j$ (when this event has positive probability), the counts $d_i=d_j=1$ have conditional covariance zero, whereas the corresponding entry of the empirical cross-product block $C_{ij}$ is $-1/(n-1)$ for $n\ge2$. Thus these blocks are not interpreted as established covariance estimators of the aggregates at random observed failure times. Instead, they provide an algebraic basis for the Delta-method constructions below, whose first-order validity is established separately. The corresponding theoretical benchmarks are provided by the joint first-order theory and the expanded-variance formulation established in \citet[Theorem~9, Remark~10, and Proposition~13]{GHTheory26}.
\end{remark}

\subsection{Full between-component correction}

Recall the between-component log-scale covariance contribution
\[
\Gamma_{\widehat A}(t)=4A.Cov\{\log\widehat S_n(t),\log\widehat H_n(t)\}.
\]
Write $H=\widehat H_n(t)$, $D_j=d_j/[n_j(n_j-d_j)]$ and
\begin{equation}
\label{eq:Q_def_main}
\widehat Q_n(t):=\sum_{j:t_j\le t}\frac{d_j}{n_j(n_j-d_j)^2}.
\end{equation}
For brevity let $Q=\widehat Q_n(t)$ and $T_2=\sum_{j:t_j\le t}D_j^2$. Work first on the regular set with $H>0$. The gradients of $\ell(x,z)=\log(1-x/z)$ and $\phi(x,z)=x/[z(z-x)]$ at the observed arguments are
\[
\alpha_i=\begin{pmatrix}-1/(n_i-d_i)\\ D_i\end{pmatrix},\qquad
\beta_j=\begin{pmatrix}1/(n_j-d_j)^2\\-d_j(2n_j-d_j)/[n_j^2(n_j-d_j)^2]\end{pmatrix}.
\]
Since $\log\widehat S_n=\sum_i\ell(d_i,n_i)$ and $\nabla_{(d_j,n_j)}\log H=\beta_j/H$ in the formal finite-dimensional calculation, the empirical blocks $C_{ij}$ motivate the mixed Delta candidate
\[
K_{ij}^{(A)}:=\alpha_i^\top C_{ij}\beta_j,\qquad
\widehat\Gamma_{\widehat A,\Delta}(t):=\frac4H\sum_{i,j:t_i,t_j\le t}K_{ij}^{(A)}.
\]
The three index cases are (see the algebraic derivation in the Appendix, Section~\ref{app:GammaA_algebra})
\[
K_{ij}^{(A)}=\begin{cases}
-\dfrac{n}{n-1}\dfrac{d_i}{n_i(n_i-d_i)^2},&i=j,\\[0.6ex]
-\dfrac{n}{n-1}D_iD_j,&i<j,\\[0.6ex]
0,&i>j.
\end{cases}
\]
Thus, using $2\sum_{i<j}D_iD_j=H^2-T_2$, we obtain the exact sample identity
\begin{equation}
\label{eq:GammaA_raw}
\widehat\Gamma_{\widehat A,\Delta}(t)
=\frac{n}{n-1}\left(-\frac{4Q}{H}-2H+\frac{2T_2}{H}\right).
\end{equation}
We use the compact first-order-equivalent definition
\begin{equation}
\label{eq:GammaA_full_main}
\boxed{
\widehat\Gamma_{\widehat A,\mathrm{full}}(t)
=-4\widehat H_n(t)^{-1}\widehat Q_n(t)-2\widehat H_n(t).
}
\end{equation}
Both candidates are set to zero whenever $H=0$ and on the exceptional samples specified in Subsection~\ref{subsec:global_well_definedness}. On the regular set with $H>0$, their exact difference is
\[
\widehat\Gamma_{\widehat A,\Delta}(t)-\widehat\Gamma_{\widehat A,\mathrm{full}}(t)
=\frac{\widehat\Gamma_{\widehat A,\mathrm{full}}(t)}{n-1}
+\frac{2nT_2}{(n-1)H}.
\]
Uniform risk-set positivity with at most $n$ summands gives $T_2=O_p(n^{-3})$. By the direct consistency of $nH$ and the plug-in limit for $n^2Q$ (Appendix Section~\ref{app:GammaA_algebra} and Lemma~\ref{lem:Q_limit}),
\[
nH\xrightarrow{P}h(t)>0,\qquad n^2Q\xrightarrow{P}J(t):=\int_0^t\frac{\lambda(u)}{y(u)^2}\,du.
\]
Hence the inverse of $H$ on the positive set is $O_p(n)$, and the difference is $O_p(n^{-2})$. On the original Greenwood scale, multiplication by $\widehat G_n(t)^2=O_p(n^{-2})$ makes this difference $O_p(n^{-4})$. In particular, the two candidates are not asserted to be exactly equal. From \eqref{eq:GammaA_full_main}, and therefore also for the raw candidate \eqref{eq:GammaA_raw},
\[
n\widehat\Gamma_{\widehat A,\mathrm{full}}(t)\xrightarrow{P}-\frac{4J(t)}{h(t)}-2h(t).
\]
We now compare the empirical limit with the joint CLT of the accompanying theory paper. Theorem~9 and Remark~10 of \citet{GHTheory26} give $c_{SH}(t)=-S(t)\{J(t)+h(t)^2/2\}$. The log-Delta method, with $\log\widehat H_n=\log\widetilde H_n-\log n$ on the positive set, gives the theoretical contribution
\begin{equation}
\label{eq:GammaA_theoretical_main}
\Gamma_{\widehat A}(t)=\frac{4c_{SH}(t)}{nS(t)h(t)}=-\frac{4J(t)}{nh(t)}-\frac{2h(t)}{n}.
\end{equation}
This agrees with the independently obtained limit of the empirical candidate.

\subsection{Full Greenwood-component correction}
\label{subsec:GammaB_construction}

The base construction in Subsection~\ref{4.3} retains the
diagonal Greenwood contribution while leaving the remaining first-order variance of the Greenwood component unaccounted for. We now
construct the corresponding count-based correction directly as
an algebraic candidate motivated by first-order bivariate Delta-method covariance propagation. Its
first-order limit is established in Appendix Lemma~\ref{lem:GammaB_limit}
and stated explicitly below. The resulting correction enters the
complete count-based estimator in Subsection~\ref{subsec:Rcorr_validation},
where Theorem~\ref{thm:Rcorr_validation} establishes agreement with the
independent first-order benchmark from the theory paper.

Write the individual Greenwood summands as
\[
\widehat H_n(t)=\sum_{j:t_j\le t}D_j,
\qquad
D_j=\phi(d_j,n_j),
\qquad
\phi(x,z):=\frac{x}{z(z-x)}.
\]
Both the observed event count \(d_j\) and the predictable risk-set size
\(n_j\) enter this two-variable map. Its gradient at the observed values is
\[
\nabla\phi(d_j,n_j)
=\begin{pmatrix}a_j\\ b_j\end{pmatrix},
\qquad
a_j:=\frac{1}{(n_j-d_j)^2},
\quad
b_j:=-\frac{d_j(2n_j-d_j)}{n_j^2(n_j-d_j)^2}.
\]
Let $\beta_j:=\nabla\phi(d_j,n_j)=(a_j,b_j)^\top$. For two distinct event times $t_i<t_j$, define the algebraic pairwise Delta candidate using the empirical cross-product block introduced above,
\[
K_{ij}^{(B)}:=\beta_i^\top C_{ij}\beta_j.
\]
Substitution of the block entries, followed by elementary simplification, gives the entirely count-based pair contribution
\[
K_{ij}^{(B)}
=
\frac{d_i d_j\{n_i^2+(2n_i-d_i)(n-n_i)\}}
{(n-1)n_i^2(n_i-d_i)^2 n_j(n_j-d_j)},
\qquad i<j.
\]
The expression is a Delta-method candidate based on the observed event counts and risk-set sizes, with aggregate first-order validity established independently in Appendix Lemma~\ref{lem:GammaB_limit}.

Under the continuous-time no-ties model, \(d_i=d_j=1\). Summing twice
the corresponding pair contributions in the count-based pair expression
gives the additional count-based correction
\begin{equation}
\label{eq:GammaB_nt_main}
\boxed{
\widehat\Gamma_{\widehat B,\mathrm{nt}}(t)
=
2
\sum_{\substack{i<j\\ t_i,t_j\le t}}
\frac{
2n n_i-n_i^2-n+n_i
}{
(n-1)n_i^2 n_j (n_i-1)^2(n_j-1)
}.
}
\end{equation}
Its first-order limit is proved in Appendix
Lemma~\ref{lem:GammaB_limit} and used in the validation in
Subsection~\ref{subsec:Rcorr_validation}.

\subsection{Corrected variance estimator and first-order validation}
\label{subsec:Rcorr_validation}

The two corrections enter the original Greenwood scale differently.
Because \(\Gamma_{\widehat A}(t)\) is a log-scale covariance contribution,
its plug-in version is multiplied by
\[
\widehat G_n(t)^2
=
\widehat S_n(t)^4\widehat H_n(t)^2.
\]
By contrast,
\(\widehat\Gamma_{\widehat B,\mathrm{nt}}(t)\) is already on the
original \(\widehat H_n\)-variance scale. On the regular set with $\widehat H_n(t)>0$, the factor
\(\widehat H_n(t)^{-2}\) associated with the logarithmic transformation
and the factor \(\widehat H_n(t)^2\) from the final back-transformation
therefore cancel, leaving the factor \(\widehat S_n(t)^4\). On the remaining samples the empirical contributions are zero under the common extension convention.

The main result of the paper is as follows.

\begin{theorem}[Corrected Greenwood variance formula]
\label{thm:Rcorr_validation}
Assume Assumptions~\ref{ass:iid-independent-censoring}--%
\ref{ass:ties}. Also assume the log-positivity condition
\eqref{eq:h_positive_log} at the fixed time point \(t\le\tau\).
Then the corrected original-scale variance estimator is
\begin{equation}
\label{eq:Rcorr_base_gammaA_gammaB}
\boxed{
\widehat R_{\mathrm{corr}}(t)
=
\widehat R_{\mathrm{base}}(t)
+
\widehat S_n(t)^4\widehat H_n(t)^2
\widehat\Gamma_{\widehat A,\mathrm{full}}(t)
+
\widehat S_n(t)^4
\widehat\Gamma_{\widehat B,\mathrm{nt}}(t).
}
\end{equation}
Equivalently,
\begin{equation}
\label{eq:Rcorr_expanded}
\boxed{
\widehat R_{\mathrm{corr}}(t)
=
\widehat S_n(t)^4
\left[
2\widehat H_n(t)^3
+
\widehat C_{\mathrm{diag}}(t)
-
4\widehat H_n(t)\widehat Q_n(t)
+
\widehat\Gamma_{\widehat B,\mathrm{nt}}(t)
\right].
}
\end{equation}
Moreover,
\[
n^3\widehat R_{\mathrm{corr}}(t)
\xrightarrow{P}
v_G(t),
\]
where
\begin{equation}
\label{eq:vG_full_main}
v_G(t)
=
S(t)^4
\left\{
v_H(t)+2h(t)^3-4h(t)J(t)
\right\}.
\end{equation}
The quantity \(v_G(t)\) is exactly the limiting CLT variance of
\[
\sqrt n\{
\widetilde G_n(t)-S(t)^2h(t)
\}.
\]
\end{theorem}

\begin{proof}
We use Theorems~5 and~9 of \citet{GHTheory26} (in particular,
Kaplan--Meier consistency from the joint CLT), and Appendix
Lemmas~\ref{lem:Q_limit}, \ref{lem:diag_limit}, and~\ref{lem:GammaB_limit}:
\[
\begin{gathered}
\widehat S_n(t)\xrightarrow{P}S(t),\qquad
n\widehat H_n(t)\xrightarrow{P}h(t),\\
n^2\widehat Q_n(t)\xrightarrow{P}J(t),\qquad
n^3\widehat C_{\mathrm{diag}}(t)\xrightarrow{P}L(t).
\end{gathered}
\]
In particular, Appendix Lemma~\ref{lem:GammaB_limit} gives the empirically derived no-ties pair limit
\begin{equation}
\label{eq:GammaB_KB_limit}
n^3\widehat\Gamma_{\widehat B,\mathrm{nt}}(t)
\xrightarrow{P}K_B(t),
\end{equation}
where $K_B(t)$ is the integral in \eqref{eq:GammaB_K_integral}. Together with the diagonal limit, this yields
\begin{equation}
\label{eq:GammaB_limit_main}
n^3\{\widehat C_{\mathrm{diag}}(t)
+\widehat\Gamma_{\widehat B,\mathrm{nt}}(t)\}
\xrightarrow{P}L(t)+K_B(t).
\end{equation}

From \eqref{eq:Rcorr_expanded},
\[
\begin{aligned}
n^3\widehat R_{\mathrm{corr}}(t)
=\widehat S_n(t)^4\big[&2\{n\widehat H_n(t)\}^3
+n^3\widehat C_{\mathrm{diag}}(t)\\
&-4\{n\widehat H_n(t)\}\{n^2\widehat Q_n(t)\}
+n^3\widehat\Gamma_{\widehat B,\mathrm{nt}}(t)\big].
\end{aligned}
\]
Consequently, by Slutsky's lemma and the continuous mapping theorem,
\[
\begin{aligned}
n^3\widehat R_{\mathrm{corr}}(t)
&\xrightarrow{P}S(t)^4\{2h(t)^3+L(t)-4h(t)J(t)+K_B(t)\}.
\end{aligned}
\]
We now identify the empirical limits with the theory-paper variances. The calculation in Appendix Lemma~\ref{lem:GammaB_limit} identifies $L(t)+K_B(t)=v_H(t)$ using Proposition~13 of \citet{GHTheory26}. Hence
\[
\begin{aligned}
S(t)^4\{2h(t)^3+L(t)-4h(t)J(t)+K_B(t)\}
&=S(t)^4\{v_H(t)+2h(t)^3-4h(t)J(t)\}\\
&=v_G(t).
\end{aligned}
\]
The last equality is the already proved theory-paper variance identity
(Corollary~11 and Remark~12 of \citet{GHTheory26}), used here as the
first-order benchmark. Thus the equality of limits completes the proof
without repeating the theory-paper influence-function calculation.
\end{proof}

For direct numerical implementation, the corrected estimator can be written
entirely in terms of observed risk-set quantities. First,
\begin{equation}
\label{eq:Cdiag_Q_explicit}
\widehat C_{\mathrm{diag}}(t)
=
\sum_{j:t_j\le t}
\frac{d_j}{n_j(n_j-d_j)^3},
\qquad
\widehat Q_n(t)
=
\sum_{j:t_j\le t}
\frac{d_j}{n_j(n_j-d_j)^2}.
\end{equation}
Under the continuous-time no-ties model,
\[
d_j=1
\]
at every observed event time, so that
\begin{equation}
\label{eq:H_C_Q_noties_explicit}
\widehat H_n(t)
=
\sum_{j:t_j\le t}\frac{1}{n_j(n_j-1)},
\;
\widehat C_{\mathrm{diag}}(t)
=
\sum_{j:t_j\le t}\frac{1}{n_j(n_j-1)^3},
\;
\widehat Q_n(t)
=
\sum_{j:t_j\le t}\frac{1}{n_j(n_j-1)^2}.
\end{equation}
Substitution of these expressions and
\eqref{eq:GammaB_nt_main} into \eqref{eq:Rcorr_expanded} gives the fully
explicit no-ties implementation formula
\begin{equation}
\label{eq:Rcorr_fully_explicit_noties}
\boxed{
\begin{aligned}
\widehat R_{\mathrm{corr}}(t)
=
\widehat S_n(t)^4
\Bigg[
&2\left\{
\sum_{j:t_j\le t}
\frac{1}{n_j(n_j-1)}
\right\}^3
+
\sum_{j:t_j\le t}
\frac{1}{n_j(n_j-1)^3}
\\[0.5ex]
&-4\left\{
\sum_{j:t_j\le t}
\frac{1}{n_j(n_j-1)}
\right\}
\left\{
\sum_{j:t_j\le t}
\frac{1}{n_j(n_j-1)^2}
\right\}
\\[0.5ex]
&+
2\sum_{\substack{i<j\\ t_i,t_j\le t}}
\frac{
2n n_i-n_i^2-n+n_i
}{
(n-1)n_i^2n_j(n_i-1)^2(n_j-1)
}
\Bigg].
\end{aligned}
}
\end{equation}
Equation~\eqref{eq:Rcorr_fully_explicit_noties} contains only directly
observable quantities: the total sample size \(n\), the risk-set sizes
\(n_j\), the observed event times through the summation ranges, and the
Kaplan--Meier factor \(\widehat S_n(t)\). It is therefore the preferred form
for direct numerical implementation under the no-ties framework, with the zero-extension convention in Subsection~\ref{subsec:global_well_definedness}.

\begin{remark}[No-censoring check]
\label{rem:no_censoring_Rcorr}
If there is no censoring, then on the regular event
\[
\mathcal D_n(t)=\{N_n(t)<n\}=\{\widehat S_n(t)>0\},
\]
we have, for $n\ge2$,
\[
\widetilde H_n(t)
=
\frac{1}{\widehat S_n(t)}-1
\]
and therefore
\[
\widetilde G_n(t)
=
\widehat S_n(t)^2\widetilde H_n(t)
=
\widehat S_n(t)\{1-\widehat S_n(t)\}.
\]
Since $S(t)>0$, the complementary event has probability $\{1-S(t)\}^n\to0$, so the fixed extension of the Greenwood component does not affect the following first-order calculation. The product identity for $\widetilde G_n(t)$ also holds on $\mathcal D_n(t)^c$ under the zero extension.
Since, in this case,
\[
\sqrt n\{\widehat S_n(t)-S(t)\}
\xRightarrow{d}
N\{0,S(t)(1-S(t))\},
\]
the ordinary Delta method applied to
\(x\mapsto x(1-x)\) gives the direct variance check
\[
v_G(t)
=
(1-2S(t))^2S(t)\{1-S(t)\}.
\]
On the other hand,
\[
h(t)=\frac{1-S(t)}{S(t)},
\qquad
J(t)=\frac12\{S(t)^{-2}-1\},
\qquad
v_H(t)=\frac{1-S(t)}{S(t)^3}.
\]
Substitution into \eqref{eq:vG_full_main} gives exactly
\[
S(t)^4
\{v_H(t)+2h(t)^3-4h(t)J(t)\}
=
(1-2S(t))^2S(t)\{1-S(t)\}.
\]
Thus the corrected formula agrees with the direct no-censoring
Delta-method calculation. In particular, $S(t)=1/2$ gives $h(t)=1>0$ but $v_G(t)=0$.
\end{remark}

\begin{corollary}[Asymptotic non-negativity]
\label{cor:Rcorr_asymptotic_nonnegative}
Under the assumptions of Theorem~\ref{thm:Rcorr_validation},
\(v_G(t)\ge0\). If \(v_G(t)>0\), then
\[
P\{\widehat R_{\mathrm{corr}}(t)\ge0\}
\longrightarrow1.
\]
No exact finite-sample non-negativity claim is made for the componentwise
plug-in formula \eqref{eq:Rcorr_expanded}.
\end{corollary}

\begin{proof}
The quantity \(v_G(t)\) is a variance and is therefore nonnegative. If it is
strictly positive, convergence in probability of
\(n^3\widehat R_{\mathrm{corr}}(t)\) to \(v_G(t)\) implies that the corrected
estimator is nonnegative with probability tending to one.
\end{proof}

\subsection{Consistency for the Actual Sampling Variance}
\label{subsec:actual_sampling_variance}

The first-order validation above is a convergence-in-probability result for the corrected estimator. To identify its target with the leading term of the actual sampling variance, we establish the following additional moment result. Here $v_G(t)$ denotes the limiting CLT variance in Corollary~11 of \citet{GHTheory26}, which agrees with \eqref{eq:vG_full_main} under the assumptions of Theorem~\ref{thm:Rcorr_validation}.

\begin{proposition}[Moment convergence]
\label{prop:G_moment_convergence}
Under Assumptions~\ref{ass:iid-independent-censoring}--\ref{ass:ties} and the global convention of Subsection~\ref{subsec:global_well_definedness}, let $t\le\tau$ be fixed and set
\[
W_n(t)=\sqrt n\{\widetilde G_n(t)-S(t)^2h(t)\}.
\]
Then $\sup_{n\ge1}E|W_n(t)|^4<\infty$, and
\[
EW_n(t)\to0,\qquad EW_n(t)^2\to v_G(t),\qquad
n^3\operatorname{Var}\{\widehat G_n(t)\}\to v_G(t).
\]
Neither $h(t)>0$ nor $v_G(t)>0$ is required.
\end{proposition}
The proof is given in Subsection~\ref{subsec:G_moment_proof}.

\begin{corollary}[Consistency for the actual sampling variance]
\label{cor:Rcorr_actual_variance}
Under the assumptions of Theorem~\ref{thm:Rcorr_validation},
\[
n^3[\widehat R_{\mathrm{corr}}(t)-\operatorname{Var}\{\widehat G_n(t)\}]
\xrightarrow{P}0.
\]
If, in addition, $v_G(t)>0$, then the actual sampling variance is strictly positive for all sufficiently large $n$, and
\[
\frac{\widehat R_{\mathrm{corr}}(t)}{\operatorname{Var}\{\widehat G_n(t)\}}
\xrightarrow{P}1.
\]
\end{corollary}
\begin{proof}
By Theorem~\ref{thm:Rcorr_validation} and Proposition~\ref{prop:G_moment_convergence}, both $n^3\widehat R_{\mathrm{corr}}(t)$ and the deterministic sequence $n^3\operatorname{Var}\{\widehat G_n(t)\}$ converge to $v_G(t)$, the former in probability. Their difference therefore converges in probability to zero. If $v_G(t)>0$, the scaled sampling variance is at least $v_G(t)/2$ for all sufficiently large $n$, and Slutsky's theorem gives the ratio convergence (with arbitrary extension at any earlier zero denominators).
\end{proof}

If $h(t)=0$, positivity of $y$ and nonnegativity of $\lambda$ imply $\lambda=0$ almost everywhere on $[0,t]$. Since $E N_i(t)=\int_0^t y(u)\lambda(u)\,du=0$, there are almost surely no observed failures up to $t$. Both $\widehat G_n(t)$ and $\widehat R_{\mathrm{corr}}(t)$ are then zero under the global convention, and the scaled difference is identically zero almost surely; no ratio is asserted. For $h(t)>0$ and $v_G(t)=0$, the scaled difference still converges to zero in probability, but neither relative consistency nor asymptotic non-negativity follows from this convergence alone.

\section{Discussion}
\label{sec:discussion}

This article develops an explicit estimator of the original-scale first-order variance approximation for Greenwood's variance estimator. Starting from its product representation,
the construction separates a directly computable diagonal base contribution
from two additional first-order terms: the between-component covariance
correction and the pairwise correction for the remaining first-order variance of the Greenwood component. Combining these contributions yields the corrected estimator
\(\widehat R_{\mathrm{corr}}(t)\), expressed entirely in observed event counts
and risk-set sizes.

The present paper and the accompanying theory paper \citet{GHTheory26}
address complementary questions. The theory paper establishes the asymptotic
linearity of the scaled Greenwood component and its joint first-order behavior
with the Kaplan--Meier estimator, thereby providing the limiting CLT variance
\(v_G(t)\). The present paper constructs both corrections from empirical cross-product blocks, simplifies the between-component candidate by first-order equivalence, and proves that its
\(n^3\)-scaled limit agrees with this independent influence-function benchmark
(Theorem~\ref{thm:Rcorr_validation}).

The validation is pointwise at a fixed time $t$ under the regularity conditions stated in Section~\ref{definitions}, including the continuous-time no-ties model and the log-positivity condition $h(t)>0$.

It establishes $n^3\widehat R_{\mathrm{corr}}(t)\xrightarrow{P}v_G(t)$, where $v_G(t)$ is the limiting CLT variance of the scaled, centered Greenwood error. Thus, the corrected estimator consistently estimates the original-scale first-order variance approximation $v_G(t)/n^3$ on the $n^3$-scale.

The first-order validation alone does not establish convergence of $n^3\operatorname{Var}\{\widehat G_n(t)\}$ to $v_G(t)$. Proposition~\ref{prop:G_moment_convergence} establishes this convergence separately, and Corollary~\ref{cor:Rcorr_actual_variance} then gives consistency for the actual sampling variance on the $n^3$-scale, as well as relative consistency if $v_G(t)>0$.

Moreover, the corrected estimator is asymptotically nonnegative whenever the limiting variance is strictly positive, whereas no exact finite-sample non-negativity property is claimed. If $v_G(t)=0$, the scaled difference convergence alone gives neither relative consistency nor a corresponding non-negativity probability tending to one.

Several methodological extensions suggest themselves. First,
the full corrected estimator and its asymptotic validation could be extended
beyond the continuous-time no-ties setting to allow tied event times. Second,
the present results are pointwise at a fixed time \(t\). A natural theoretical
extension would be to develop uniform or process-level results over a time
interval, thereby moving from pointwise first-order asymptotics to a functional
treatment of the sampling variability of Greenwood's variance estimator.

A complementary methodological direction would be to investigate an
alternative count-based estimator obtained by directly plugging empirical
risk-set quantities and Nelson--Aalen hazard increments into the expanded
first-order variance formula derived in the accompanying theory paper.
Such an approach would start from the theoretical limiting variance rather
than from the local Delta-method covariance construction pursued here.
Of particular interest would be the algebraic relationship between the
resulting plug-in estimator and the present corrected formula, as well as a
comparison of their finite-sample performance, particularly in terms of
bias, sampling variability, and mean squared error when estimating the
sampling variance of Greenwood's estimator. Further questions concern their
higher-order asymptotic behavior and nonnegativity properties.

A broader inferential direction concerns variability in
individual time-to-event outcomes themselves. The present article studies the
sampling variability of Greenwood's variance estimator and therefore the
uncertainty associated with estimation of a pointwise survival probability. A
different question is the variability of survival times in the underlying
population. Addressing that question requires variability or dispersion measures
defined directly as functionals of the survival-time distribution rather than
as sampling variances of estimators. Developing suitable estimators and
corresponding inferential procedures for such functionals under censoring would
therefore constitute a natural further direction. Such methods could ultimately
provide a basis for comparing outcome variability between treatment groups in
time-to-event settings and for investigating heterogeneity beyond differences
in mean or pointwise survival.

\section{Appendix}
\label{auxFundamental2}
\label{auxTech}

This appendix records the auxiliary first-order plug-in limits used in
Sections~\ref{proof} and~\ref{sec:improvement_base_formula}. The full regularity
of the scaled Greenwood component is established separately in \citet{GHTheory26}. Subsection~\ref{subsec:G_moment_proof} provides the additional moment-convergence argument for the Greenwood estimator itself.

\begin{lemma}[Greenwood plug-in for the log-Kaplan--Meier variance]
\label{lem:greenwood_log_variance_base}
Under the standing fixed-time no-ties assumptions,
\[
A.Var\{\log\widehat S_n(t)\}
=
\frac{h(t)}{n},
\qquad
h(t)
=
\int_0^t\frac{\lambda(u)}{y(u)}\,du.
\]
Moreover,
\[
n\widehat H_n(t)\xrightarrow{P}h(t),
\]
and therefore
\[
n\left[
\widehat H_n(t)-A.Var\{\log\widehat S_n(t)\}
\right]
\xrightarrow{P}0.
\]
Thus \(\widehat H_n(t)\) is a sample-based plug-in estimator of the
theoretical original-scale first-order variance
\(A.Var\{\log\widehat S_n(t)\}=h(t)/n\).
\end{lemma}

\begin{proof}
The standard first-order theory of the Kaplan--Meier estimator (see the first component of Theorem~9 of \citet{GHTheory26}) gives
\[
\sqrt n\{\widehat S_n(t)-S(t)\}
\xRightarrow{d}
N\{0,S(t)^2h(t)\}.
\]
Since \(S(t)>0\) under the standing assumptions, the Delta method applied to
\(x\mapsto\log x\) yields
\[
\sqrt n\{\log\widehat S_n(t)-\log S(t)\}
\xRightarrow{d}
N\{0,h(t)\}.
\]
By the original-scale convention of Section~\ref{subsec:Avar_convention}, this is exactly
\[
A.Var\{\log\widehat S_n(t)\}=\frac{h(t)}{n}.
\]
The convergence
\[
n\widehat H_n(t)\xrightarrow{P}h(t)
\]
is already contained in \citet[Theorem~4.2.2]{Gill1980} and is also proved
under the present assumptions in \citet[Theorem~5]{GHTheory26}.
Combining the two relations gives
\[
n\left[
\widehat H_n(t)-A.Var\{\log\widehat S_n(t)\}
\right]
=
n\widehat H_n(t)-h(t)
\xrightarrow{P}0,
\]
which proves the stated plug-in property.
\end{proof}

\begin{lemma}[Limit of the diagonal Greenwood plug-in]
\label{lem:diag_limit}
Under the standing no-ties assumptions, define
\[
\widehat C_{\mathrm{diag}}(t)
=
\sum_{j:t_j\le t}
\frac{1}{n_j(n_j-1)^3}.
\]
Then
\[
n^3\widehat C_{\mathrm{diag}}(t)
\xrightarrow{P}
L(t)
=
\int_0^t\frac{\lambda(u)}{y(u)^3}\,du.
\]
\end{lemma}

\begin{proof}
By the uniform law of large numbers for the empirical at-risk process and the
positivity assumption, there exists \(c>0\) such that
\[
\mathcal E_n(t)
:=
\left\{
\inf_{0\le u\le t}\frac{Y_n(u)}{n}\ge c
\right\}
\]
has probability tending to one. For all sufficiently large \(n\),
\(\mathcal E_n(t)\subseteq\mathcal D_n(t)\). Hence, on
\(\mathcal E_n(t)\), the no-ties counting-process identity gives
\[
\widehat C_{\mathrm{diag}}(t)
=
\int_0^t
\frac{1}{Y_n(u)\{Y_n(u)-1\}^3}\,dN_n(u).
\]
Therefore, on this event,
\[
n^3\widehat C_{\mathrm{diag}}(t)
=
\int_0^t a_{n,3}(u)\,\frac{dN_n(u)}{n},
\qquad
a_{n,3}(u)
:=
\begin{cases}
\dfrac{n^4}{Y_n(u)\{Y_n(u)-1\}^3},&Y_n(u)\ge2,\\
0,&Y_n(u)<2.
\end{cases}
\]
Uniform convergence of \(Y_n/n\) to \(y\), together with uniform positivity,
implies
\[
\sup_{0\le u\le t}
\left|
a_{n,3}(u)-\frac{1}{y(u)^4}
\right|
\xrightarrow{P}0.
\]
Since \(N_n(t)/n\le1\), replacing \(a_{n,3}\) by \(y^{-4}\) in the
preceding integral changes it by \(o_p(1)\). Moreover,
\[
\int_0^t\frac{1}{y(u)^4}\,\frac{dN_n(u)}{n}
=
\frac1n\sum_{i=1}^n
\int_0^t\frac{1}{y(u)^4}\,dN_i(u).
\]
The summands are i.i.d. and bounded, and the multiplicative-intensity
representation gives
\[
E\left\{
\int_0^t\frac{1}{y(u)^4}\,dN_i(u)
\right\}
=
\int_0^t\frac{1}{y(u)^4}y(u)\lambda(u)\,du
=
L(t).
\]
The ordinary law of large numbers therefore yields the asserted convergence.
Since the probability of \(\mathcal E_n(t)\) tends to one, the arbitrary
extension on \(\mathcal D_n(t)^c\) is asymptotically irrelevant.
\end{proof}

\begin{lemma}[Limit of the between-component event-count plug-in]
\label{lem:Q_limit}
Under the standing fixed-time continuous-time no-ties assumptions, define
\[
\widehat Q_n(t)
=
\sum_{j:t_j\le t}
\frac{1}{n_j(n_j-1)^2}.
\]
Then
\[
n^2\widehat Q_n(t)
\xrightarrow{P}
J(t)
=
\int_0^t\frac{\lambda(u)}{y(u)^2}\,du.
\]
\end{lemma}

\begin{proof}
Let \(\mathcal E_n(t)\) be the high-probability event used in the proof of
Lemma~\ref{lem:diag_limit}. For all sufficiently large \(n\), on
\(\mathcal E_n(t)\subseteq\mathcal D_n(t)\),
\[
\widehat Q_n(t)
=
\int_0^t
\frac{1}{Y_n(u)\{Y_n(u)-1\}^2}\,dN_n(u).
\]
Consequently,
\[
n^2\widehat Q_n(t)
=
\int_0^t a_{n,2}(u)\,\frac{dN_n(u)}{n},
\qquad
a_{n,2}(u)
:=
\begin{cases}
\dfrac{n^3}{Y_n(u)\{Y_n(u)-1\}^2},&Y_n(u)\ge2,\\
0,&Y_n(u)<2.
\end{cases}
\]
Uniform convergence and positivity imply
\[
\sup_{0\le u\le t}
\left|
a_{n,2}(u)-\frac{1}{y(u)^3}
\right|
\xrightarrow{P}0.
\]
Again using \(N_n(t)/n\le1\), the random integrand may be replaced by
\(y^{-3}\) with an \(o_p(1)\) error. Finally,
\[
\int_0^t\frac{1}{y(u)^3}\,\frac{dN_n(u)}{n}
=
\frac1n\sum_{i=1}^n
\int_0^t\frac{1}{y(u)^3}\,dN_i(u)
\xrightarrow{P}
\int_0^t\frac{1}{y(u)^3}y(u)\lambda(u)\,du
=
J(t)
\]
by the ordinary law of large numbers and the multiplicative-intensity
representation. Since the probability of \(\mathcal E_n(t)\) tends to one,
the arbitrary extension outside \(\mathcal D_n(t)\) does not affect the limit.
\end{proof}

\begin{lemma}[First-order limit of the no-ties pair correction]
\label{lem:GammaB_limit}
Under the standing fixed-time continuous-time no-ties assumptions, let
\(\widehat\Gamma_{\widehat B,\mathrm{nt}}(t)\) be defined by
\eqref{eq:GammaB_nt_main}. Then
\[
n^3\widehat\Gamma_{\widehat B,\mathrm{nt}}(t)
\xrightarrow{P}
K_B(t):=2\int_0^t\int_0^r
\frac{\lambda(s)\lambda(r)\{2-y(s)\}}{y(s)^2y(r)}\,ds\,dr.
\]
Consequently, by Lemma~\ref{lem:diag_limit},
\[
n^3\{\widehat C_{\mathrm{diag}}(t)
+\widehat\Gamma_{\widehat B,\mathrm{nt}}(t)\}
\xrightarrow{P}L(t)+K_B(t).
\]
Moreover, the expanded-variance result of the theory paper identifies $L(t)+K_B(t)=v_H(t)$.
\end{lemma}

\begin{proof}
For all sufficiently large $n$, we work on the high-probability event \(\mathcal E_n(t)\) introduced in the proof of Lemma~\ref{lem:diag_limit}, on which all risk sets under consideration satisfy \(n_j/n\ge c>0\). Write
\[
y_n(u):=\frac{Y_n(u)}{n},\qquad
\mu_n(du):=\frac{1}{n}\,dN_n(u),\qquad
\mu(du):=y(u)\lambda(u)\,du.
\]
At an observed event time \(t_j\), let \(x_j=n_j/n=y_n(t_j-)\).
The summand in \eqref{eq:GammaB_nt_main} has the exact representation
\[
\frac{2nn_i-n_i^2-n+n_i}
{(n-1)n_i^2 n_j(n_i-1)^2(n_j-1)}
=\frac{1}{n^5}\,k_n(x_i,x_j),
\]
where
\[
k_n(x,z):=
\frac{2x-x^2+(x-1)/n}
{(1-1/n)x^2z(x-1/n)^2(z-1/n)}.
\]
Uniform convergence of \(y_n\) to \(y\), positivity of \(y\), and the
bounds defining \(\mathcal E_n(t)\) imply, for the following nonnegative pair error (defined as zero on $\mathcal E_n(t)^c$, when $N_n(t)<2$, or when $n<2/c$),
\[
\sup_{\substack{i<j\\t_i,t_j\le t}}
\left|k_n(x_i,x_j)-k\{y(t_i),y(t_j)\}\right|
\xrightarrow{P}0,
\qquad
k(x,z):=\frac{2x-x^2}{x^4z^2}
=\frac{2-x}{x^3z^2}.
\]
If \(y\) has jumps due to censoring atoms, the pointwise values at
these jumps are immaterial: observed failure times hit any fixed
censoring atom with probability zero, and the limiting event measure
\(\mu\) has no atoms.
Because \(N_n(t)/n\le 1\), replacement of \(k_n\) by \(k\) in the
normalized double sum incurs an \(o_p(1)\) error. Thus
\begin{equation}
\label{eq:GammaB_pair_asymptotic}
n^3\widehat\Gamma_{\widehat B,\mathrm{nt}}(t)
=
\frac{2}{n^2}
\sum_{\substack{i<j\\t_i,t_j\le t}}
k\{y(t_i),y(t_j)\}+o_p(1)
=2\!\iint_{0\le s<r\le t}
k\{y(s),y(r)\}\,\mu_n(ds)\mu_n(dr)+o_p(1).
\end{equation}
The empirical event measure \(\mu_n\) converges weakly in probability
 to \(\mu\) by the ordinary law of large numbers for the i.i.d.
 observed-event subdistribution. As these are finite measures on
 \([0,t]\), their product measures also converge weakly in probability,
 \(\mu_n\otimes\mu_n\Rightarrow\mu\otimes\mu\).
 Define the bounded function
 \[
 f(s,r):=\mathbf1_{\{0\le s<r\le t\}}\,k\{y(s),y(r)\}.
 \]
 Its discontinuities are contained in the diagonal, the domain boundaries,
 and horizontal or vertical lines corresponding to discontinuities of
 \(y\). These are \((\mu\otimes\mu)\)-null sets: \(y\) has at
 most countably many discontinuities, and the event measure
 \(\mu(du)=y(u)\lambda(u)\,du\) is absolutely continuous.
 The Portmanteau theorem for bounded functions that are continuous
 \((\mu\otimes\mu)\)-almost everywhere, applied to these product
 measures, therefore yields the convergence in
 \eqref{eq:GammaB_pair_asymptotic}:
\begin{equation}
\label{eq:GammaB_K_integral}
\begin{aligned}
n^3\widehat\Gamma_{\widehat B,\mathrm{nt}}(t)
&\xrightarrow{P}
2\iint_{0\le s<r\le t}
\frac{2-y(s)}{y(s)^3y(r)^2}
\,\mu(ds)\mu(dr)\\
&=2\int_0^t\int_0^r
\frac{\lambda(s)\lambda(r)\{2-y(s)\}}
{y(s)^2y(r)}\,ds\,dr.
\end{aligned}
\end{equation}
The empirical pair limit is therefore $K_B(t)$ as defined in the lemma. It remains to compare it with the theory-paper variance.
The expanded-variance result for \(\widetilde H_n\) in the separate
theory paper (Proposition~13 there) gives
\[
\begin{aligned}
v_H(t)-L(t)
={}&
\int_0^t\int_0^t
\frac{\lambda(s)\lambda(r)}{y(s)^2y(r)^2}
\{y(s\vee r)-y(s)y(r)\}\,ds\,dr\\
&+2\int_0^t\frac{\lambda(r)}{y(r)}
\left\{\int_0^r\frac{\lambda(s)}{y(s)^2}\,ds\right\}dr.
\end{aligned}
\]
The first double integral is symmetric. On \(s<r\), its bracket equals
\(y(r)\{1-y(s)\}\); splitting it into the two triangles, and
combining with the second term, gives
\[
\begin{aligned}
v_H(t)-L(t)
&=2\int_0^t\int_0^r
\frac{\lambda(s)\lambda(r)\{1-y(s)\}}
{y(s)^2y(r)}\,ds\,dr
+2\int_0^t\int_0^r
\frac{\lambda(s)\lambda(r)}{y(s)^2y(r)}\,ds\,dr\\
&=2\int_0^t\int_0^r
\frac{\lambda(s)\lambda(r)\{2-y(s)\}}
{y(s)^2y(r)}\,ds\,dr.
\end{aligned}
\]
This agrees with the independently obtained limit in \eqref{eq:GammaB_K_integral}. The combined empirical limit follows from the diagonal convergence in Lemma~\ref{lem:diag_limit}, and its identification with $v_H(t)$ follows from the displayed comparison.
\end{proof}

\subsection{Algebraic derivation of the between-component Delta candidate}
\label{app:GammaA_algebra}
With $H=\widehat H_n(t)$, $Q=\widehat Q_n(t)$ and $T_2=\sum_j D_j^2$ (sums over failure times up to $t$), on the regular set with $H>0$, the gradients and the empirical cross-product blocks from Section~\ref{sec:improvement_base_formula} satisfy
\[
\begin{aligned}
\alpha_i^\top w_i&=0,&\qquad \beta_j^\top w_j&=-D_j,\\
\alpha_i^\top(1,1)^\top&=-1/n_i,&
\alpha_i^\top(-d_i,n-n_i)^\top&=nD_i.
\end{aligned}
\]
For $i=j$, direct multiplication of the diagonal block gives
\[
\alpha_i^\top C_{ii}=\frac{n}{n-1}(-d_i/n_i,0),
\]
and therefore $K_{ii}^{(A)}=-nd_i/\{(n-1)n_i(n_i-d_i)^2\}$.
For $i<j$, the block factorizes as $C_{ij}=(-d_i,n-n_i)^\top w_j^\top/(n-1)$, yielding $K_{ij}^{(A)}=-nD_iD_j/(n-1)$.
For $i>j$, $C_{ij}=w_i(-d_j,n-n_j)/(n-1)$, which gives $K_{ij}^{(A)}=0$.
Consequently,
\[
\frac4H\sum_{i,j}K_{ij}^{(A)}=-\frac{4n}{(n-1)H}\left(Q+\sum_{i<j}D_iD_j\right),
\]
which is precisely \eqref{eq:GammaA_raw}. The factor $4$ is the factor in the log-scale variance decomposition; no symmetry factor is added to the upper triangle.

For completeness, the consistency of $nH$ used in the empirical construction follows directly from the same argument as Lemmas~\ref{lem:diag_limit} and~\ref{lem:Q_limit}. On $\mathcal E_n(t)$ for all sufficiently large $n$,
\[
nH=\int_0^t\frac{1}{y_n(u)\{y_n(u)-1/n\}}\,\mu_n(du),
\]
with $y_n,\mu_n$ as in Lemma~\ref{lem:GammaB_limit}. By uniform risk-set convergence and positivity the integrand converges uniformly in probability on this event to $y^{-2}$. Since $\mu_n([0,t])\le1$, the replacement error is $o_p(1)$. The ordinary i.i.d. law of large numbers for the bounded deterministic weight $y^{-2}$ gives the limit $\int_0^t y(u)^{-2}\mu(du)=h(t)$. Together with Lemma~\ref{lem:Q_limit}, this is the empirical consistency used in Section~\ref{sec:improvement_base_formula}.

On $\mathcal E_n(t)$ for all sufficiently large $n$, $D_j=O(n^{-2})$ uniformly over the at most $n$ failure times, so $T_2=O_p(n^{-3})$. The positivity of $h(t)$ gives $H^{-1}=O_p(n)$ when the inverse is set to zero on $H=0$. It follows that $Q/H=O_p(n^{-1})$, $T_2/H=O_p(n^{-2})$, and the exact difference of the raw and compact candidates is $O_p(n^{-2})$. Its original-scale counterpart is $O_p(n^{-4})$, because $\widehat G_n(t)^2=O_p(n^{-2})$.

\subsection{Moment Convergence of Greenwood's Variance Estimator}
\label{subsec:G_moment_proof}

\begin{proof}[Proof of Proposition~\ref{prop:G_moment_convergence}]
Fix $t\le\tau$ and suppress time arguments where convenient. Set $G^*=S(t)^2h(t)$, $p=y(t)>0$, $c=p/2$, and $y_n(u)=Y_n(u)/n$. Define
\[
\delta_n=\sup_{0\le u\le t}|y_n(u)-y(u)|,\qquad
\mathcal E_n=\{\inf_{0\le u\le t}y_n(u)\ge c\}.
\]
All constants below are deterministic and independent of $n$. The Dvoretzky--Kiefer--Wolfowitz inequality \citep{Massart1990} applies to the empirical distribution function $F_n$ of $X_i$, including distributions with atoms. Since $y_n(u)=1-F_n(u-)$ and $y(u)=1-F(u-)$, their difference is controlled by the empirical distribution function supremum. If infinite observation times are allowed, use $X_i\wedge(t+1)$ instead, without changing the at-risk indicators up to $t$. The supremum is measurable by left-continuity, using a countable dense grid and the endpoints. Hence, for $x>0$,
\[
P(\delta_n>x)\le2e^{-2nx^2},\qquad
E\delta_n^4=4\int_0^\infty x^3P(\delta_n>x)\,dx
\le8\int_0^\infty x^3e^{-2nx^2}\,dx=n^{-2}.
\]
By monotonicity, $\mathcal E_n=\{y_n(t)\ge c\}$ and $\mathcal E_n^c\subseteq\{\delta_n>c\}$, so $P(\mathcal E_n^c)\le2e^{-2nc^2}$. Choose $n_0\ge2$ with $n_0c\ge2$. For $n\ge n_0$, on $\mathcal E_n$ every failure risk set satisfies $n_j\ge cn$ and $n_j-1\ge cn/2$. Thus $\mathcal E_n\subseteq\mathcal D_n(t)$ and $\widehat S_n(t)>0$ on $\mathcal E_n$.

For $k=1,2$ define $f_k(Z)=y(X)^{-k}$ on $\{X\le t,\Delta=1\}$ and zero elsewhere. These weights are bounded by $p^{-k}$. By independent censoring and absolute continuity of the event time, $E\{dN_i(u)\}=y(u)\lambda(u)\,du$ on $[0,t]$, including in the presence of censoring atoms (there are no failure atoms). Consequently $Ef_1(Z)=\Lambda(t)$ and $Ef_2(Z)=h(t)$. Write
\[
B_{k,n}=n^{-1}\sum_{i=1}^n\xi_{k,i},\qquad
\xi_{k,i}=f_k(Z_i)-Ef_k(Z).
\]
The summands are centered, bounded and i.i.d.; expansion of the fourth power gives
\[
E(\sqrt n B_{k,n})^4
=\frac{E\xi_{k,1}^4}{n}+\frac{3(n-1)}{n}(E\xi_{k,1}^2)^2\le C_k<\infty.
\]
No independence of $B_{k,n}$ and $\delta_n$ is needed.

For $n\ge n_0$ on $\mathcal E_n$, let $U_n=n\sum_{j:t_j\le t}n_j^{-2}$. The function $x\mapsto x^{-2}$ is Lipschitz on $[c,1]$ with constant $2c^{-3}$. Comparing the at most $n$ failure terms with the empirical mean of $f_2$ yields $|U_n-h(t)|\le2c^{-3}\delta_n+|B_{2,n}|$. Moreover,
\[
0\le\widetilde H_n-U_n
=n\sum_{j:t_j\le t}\frac{1}{n_j^2(n_j-1)}\le\frac{2}{nc^3}.
\]
It follows that $|\widetilde H_n-h(t)|\le2c^{-3}\delta_n+|B_{2,n}|+2/(nc^3)$ on $\mathcal E_n$. By $(a+b+d)^4\le27(a^4+b^4+d^4)$ for nonnegative terms and the preceding fourth-moment bounds,
\[
\sup_{n\ge n_0}E\{n^2|\widetilde H_n-h(t)|^4\mathbf1_{\mathcal E_n}\}<\infty.
\]

Similarly, for the Nelson--Aalen estimator $\widehat\Lambda_n=\sum_{j:t_j\le t}n_j^{-1}$, the Lipschitz constant $c^{-2}$ of $x\mapsto x^{-1}$ and the empirical mean of $f_1$ give $|\widehat\Lambda_n-\Lambda(t)|\le c^{-2}\delta_n+|B_{1,n}|$ on $\mathcal E_n$. For $0\le z\le1/2$, $0\le-\log(1-z)-z=\int_0^z s/(1-s)\,ds\le z^2$. Therefore
\[
0\le-\log\widehat S_n-\widehat\Lambda_n
\le\sum_{j:t_j\le t}n_j^{-2}\le\frac{1}{nc^2}
\]
on $\mathcal E_n$. Since $S(t)=e^{-\Lambda(t)}$ and $x\mapsto e^{-x}$ is Lipschitz with constant one on $[0,\infty)$,
\[
|\widehat S_n-S(t)|\le c^{-2}\delta_n+|B_{1,n}|+\frac{1}{nc^2}
\quad\text{on }\mathcal E_n.
\]
The same fourth-moment argument gives $\sup_{n\ge n_0}E\{n^2|\widehat S_n-S(t)|^4\mathbf1_{\mathcal E_n}\}<\infty$. The logarithm has been used only on $\mathcal E_n$; no global log-moment bound is required.

As $0\le\widehat S_n,S\le1$ and $h(t)\le\Lambda(t)/p<\infty$, the product decomposition gives
\[
|\widetilde G_n-G^*|\le|\widetilde H_n-h(t)|+2h(t)|\widehat S_n-S(t)|.
\]
Using $(a+b)^4\le8(a^4+b^4)$ and the preceding bounds, we obtain $\sup_{n\ge n_0}E\{|W_n|^4\mathbf1_{\mathcal E_n}\}<\infty$.

For the entire complement $\mathcal E_n^c$, a global bound is available. Under no ties, the integer risk-set sizes at successive failures decrease strictly. On $\mathcal D_n(t)$ they are distinct members of $\{2,\ldots,n\}$, so for $n\ge2$,
\[
0\le\widehat H_n\le\sum_{r=2}^n\frac{1}{r(r-1)}=1-\frac1n.
\]
By the global zero extension (including $n=1$), $0\le\widehat H_n\le1$ and $0\le\widetilde G_n\le n$ almost surely. Thus, for $n\ge n_0$,
\[
E\{|W_n|^4\mathbf1_{\mathcal E_n^c}\}
\le n^2(n+|G^*|)^4P(\mathcal E_n^c)
\le2n^2(n+|G^*|)^4e^{-2nc^2}\to0.
\]
This controls also regular samples with small risk sets. The finitely many $n<n_0$ have finite fourth moments by the same global bound. Combining both events gives $C=\sup_{n\ge1}E|W_n|^4<\infty$.

For $K>0$, $E\{W_n^2\mathbf1_{\{W_n^2>K\}}\}\le C/K$ and $E\{|W_n|\mathbf1_{\{|W_n|>K\}}\}\le C^{1/2}/K$. Thus both $W_n$ and $W_n^2$ are uniformly integrable. Corollary~11 of \citet{GHTheory26} gives $W_n\xRightarrow{d}N\{0,v_G(t)\}$, with the point mass at zero allowed. Weak convergence and uniform integrability yield $EW_n\to0$ and $EW_n^2\to v_G(t)$. Finally, since $G^*$ is deterministic,
\[
n^3\operatorname{Var}\{\widehat G_n(t)\}=\operatorname{Var}(W_n)
=EW_n^2-(EW_n)^2\to v_G(t).
\]
\end{proof}

\section*{Statements and Declarations}

\noindent\textbf{Competing interests.} The author declares no competing interests.

\noindent\textbf{Funding.} No specific funding was received for this work.

\noindent\textbf{Data availability.} No data were generated or analysed in this study.

\noindent\textbf{Use of generative AI.} The research question, the underlying methodological concept, and substantial parts of the mathematical derivations and proofs were developed by the author. Generative artificial intelligence (ChatGPT, OpenAI) was subsequently used to assist in identifying and addressing gaps in the proofs, verifying mathematical arguments, and revising and systematically reviewing the manuscript. The author assumes full responsibility for the scientific content and conclusions of this work.

\bibliographystyle{chicago}
\setlength{\bibsep}{2pt}
\bibliography{references2}

@book{Andersen1993,
  author = {Andersen, Per Kragh and Borgan, {\O}rnulf and Gill, Richard D. and Keiding, Niels},
  title = {Statistical Models Based on Counting Processes},
  year = {1993},
  publisher = {Springer},
  address = {New York}
}

@article{BreslowCrowley1974,
  author = {Breslow, Norman E. and Crowley, John},
  title = {A Large Sample Study of the Life Table and Product Limit Estimates under Random Censorship},
  journal = {The Annals of Statistics},
  year = {1974},
  volume = {2},
  number = {3},
  pages = {437--453}
}

@article{Chu20,
  author = {Chung, Wendy K. and Erion, Karel and Florez, Jose C. and Hattersley, Andrew T. and Hivert, Marie-France and Lee, Christine G. and McCarthy, Mark I. and Nolan, John J. and Norris, Jill M. and Pearson, Ewan R. and Philipson, Louis and McElvaine, Allison T. and Cefalu, William T. and Rich, Stephen S. and Franks, Paul W.},
  title = {Precision Medicine in Diabetes: A Consensus Report from the {American Diabetes Association} ({ADA}) and the {European Association for the Study of Diabetes} ({EASD})},
  journal = {Diabetes Care}, year = {2020}, volume = {43}, number = {7}, pages = {1617--1635}, doi = {10.2337/dci20-0022}
}

@article{Den20,
  author = {Dennis, J. M.},
  title = {Precision Medicine in Type 2 Diabetes: Using Individualized Prediction Models to Optimize Selection of Treatment},
  journal = {Diabetes}, year = {2020}, volume = {69}, number = {10}, pages = {2075--2085}
}

@book{Fleming1991,
  author = {Fleming, T. R. and Harrington, D. P.},
  title = {Counting Processes and Survival Analysis},
  year = {1991}, publisher = {Wiley}, address = {New York}
}

@book{Gill1980,
  author = {Gill, R. D.}, title = {Censoring and Stochastic Integrals},
  year = {1980}, series = {Mathematical Centre Tracts}, number = {124},
  publisher = {Mathematisch Centrum}, address = {Amsterdam}
}

@article{Gill1983,
  author = {Gill, R. D.},
  title = {Large Sample Behaviour of the Product-Limit Estimator on the Whole Line},
  journal = {The Annals of Statistics}, year = {1983}, volume = {11}, number = {1}, pages = {49--58}
}

@techreport{Gre26,
  author = {Greenwood, M.}, title = {The Natural Duration of Cancer},
  institution = {His Majesty's Stationery Office}, address = {London},
  type = {Reports on Public Health and Medical Subjects}, number = {33}, year = {1926}
}

@article{Kap58,
  author = {Kaplan, E. L. and Meier, P.},
  title = {Nonparametric Estimation from Incomplete Observations},
  journal = {Journal of the American Statistical Association},
  year = {1958}, volume = {53}, number = {282}, pages = {457--481}
}

@article{Kro22,
  author = {Krook, A. and Mulder, H.},
  title = {Pinpointing Precision Medicine for Diabetes Mellitus},
  journal = {Diabetologia}, year = {2022}, volume = {65}, number = {11}, pages = {1755--1757}
}

@article{Kus23,
  author = {Kuss, O. and Opitz, M. E. and Brandstetter, L. V. and Schlesinger, S. and Roden, M. and Hoyer, A.},
  title = {How Amenable Is Type 2 Diabetes Treatment for Precision Diabetology? A Meta-Regression of Glycaemic Control Data from 174 Randomised Trials},
  journal = {Diabetologia}, year = {2023}, volume = {66}, pages = {1622--1632}
}

@article{Nol22,
  author = {Nolan, J. J. and Kahkoska, A. R. and Semnani-Azad, Z. and Hivert, M.-F. and Ji, L. and Mohan, V. and Stacey, D. and Dabelea, D. and Rich, S. S. and Mather, K. J.},
  title = {{ADA/EASD} Precision Medicine in Diabetes Initiative: An International Perspective and Future Vision for Precision Medicine in Diabetes},
  journal = {Diabetes Care}, year = {2022}, volume = {45}, number = {2}, pages = {261--266}
}

@book{Rao73,
  author = {Rao, C. R.}, title = {Linear Statistical Inference and Its Applications},
  edition = {2nd}, year = {1973}, publisher = {Wiley}, address = {New York}
}

@article{Rei81,
  author = {Reid, N.}, title = {Influence Functions for Censored Data},
  journal = {The Annals of Statistics}, year = {1981}, volume = {9}, number = {1}, pages = {78--92}
}

@unpublished{GHTheory26,
  author = {Rodenkirchen, J.},
  title = {First-Order Asymptotics of {Greenwood's} Variance Estimator},
  year = {2026}, note = {Preprint, \href{https://arxiv.org/abs/2611.41204}{arXiv:2611.41204 [math.ST]}}
}

@article{Var24,
  author = {Vargas, K. G. and R{\"u}tten, T. and Siemes, B. and others},
  title = {Assessing the Potential for Precision Medicine in Body Weight Reduction with Regard to Type 2 Diabetes Mellitus Therapies: A Meta-Regression Analysis of 120 Randomized Controlled Trials},
  journal = {Diabetes, Obesity and Metabolism}, year = {2024}, volume = {26}, pages = {2139--2146}
}

@article{Massart1990,
  author = {Massart, Pascal},
  title = {The Tight Constant in the {Dvoretzky--Kiefer--Wolfowitz} Inequality},
  journal = {The Annals of Probability},
  year = {1990},
  volume = {18},
  number = {3},
  pages = {1269--1283}
}

\end{document}